\documentclass[a4paper,10pt,reqno]{amsart}

\usepackage{float}
\usepackage{enumitem}
\usepackage{stmaryrd}
\usepackage{amsmath,amssymb,amsfonts,amsthm}
\usepackage{graphicx,hyperref}
\usepackage[all,cmtip]{xy}
\usepackage{enumitem,color}
\usepackage[latin1]{inputenc}

\newcommand{\tabref}[1]{Table~\ref{#1}}
\newcommand{\propref}[1]{Proposition~\ref{#1}}

\newcommand{\lemref}[1]{Lemma~\ref{#1}}
\newcommand{\remref}[1]{Remark~\ref{#1}}
\newcommand{\parref}[1]{Paragraph~\ref{#1}}
\newcommand{\teoref}[1]{Theorem~\ref{#1}}
\newtheorem{theorem}{Theorem}[section]
\newtheorem{proposition}[theorem]{Proposition}
\newtheorem{lemma}[theorem]{Lemma}

\theoremstyle{definition}
\newtheorem{definition}[theorem]{Definition}
\newtheorem{remark}[theorem]{Remark}

\newtheorem{parrafo}[theorem]{}

\newcommand{\R}{\mathbb{R}}
\newcommand{\C}{\mathbb{C}}
\newcommand{\Z}{\mathbb{Z}}
\renewcommand{\H}{\mathbb{H}}
\newcommand{\N}{\mathbb{N}}
\newcommand{\g}{\mathfrak{g}}
\newcommand{\h}{\mathfrak{h}}
\newcommand{\n}{\mathfrak{n}}
\newcommand{\z}{\mathfrak{z}}
\newcommand{\m}{\mathfrak{m}}

\newcommand{\spp}{\mathfrak{sp}}
\newcommand{\uu}{\mathfrak{u}}
\newcommand{\so}{\mathfrak{so}}

\newcommand{\tr}{\mathrm{tr}}

\newcommand{\pr}{\mathrm{pr}}
\newcommand{\Om}{\Omega}
\newcommand{\BOm}{\underline\Omega}

\newcommand{\Iso}{\mathrm{Iso}}
\newcommand{\id}{\mathrm{id}}
\renewcommand{\Re}{\mathrm{Re}}

\newcommand{\Ker}{\mathrm{Ker}}

\newcommand{\menos}{\setminus}
\newcommand{\depth}[1]{d(#1)}

\newcommand{\ib}[2]{{#1}_{_{#2}}}

\newcommand{\lau}[4]{{#1}^{^{#2}}_{_{#3}}{\left( #4 \right)}}
\newcommand{\coho}[3]{{#1}^{^{#2}}{\left( #3 \right)}}

\newcommand{\negro}{\color{black}}

\title{Cohomology of regular spherical actions. }
\date{\today}

\author[J.I.~Royo Prieto]{Jos\'{e} Ignacio Royo Prieto}
\address{Matematika Saila\\ Zientzia eta Teknologia Fakultatea\\ University of the Basque Country UPV/EHU\\ Barrio Sarriena s/n\\ 48940 Leioa\\Spain.}
\email{joseignacio.royo@ehu.eus}
\thanks{Partially supported by Ministerio de Ciencia, Spain, grant PID2019-105621GB-I00.}

\author[M.~Saralegi Aranguren]{Martintxo Saralegi-Aranguren}
\address{F\`{e}d\`{e}ration CNRS\\  Nord-Pas-de-Calais FR 2956\\ UPRES-EA 2462 LML\\
Facult\'e Jean Perrin\\ Universit\'{e} d'Artois\\   Rue Jean Souvraz SP 18\\ 62 307 Lens Cedex, France.}
\email{martin.saraleguiaranguren@univ-artois.fr}

\author[R.~Wolak]{Robert Wolak}
\address{Instytut Matematyki\\ Uniwersytet Jagiellonski\\
	ul. prof. Stanis{\l}awa {\L}ojasiewicza 6
	30-348 Krak\'ow,  Poland.}
\email{robert.wolak@im.uj.edu.pl}

\begin{document}

\thispagestyle{empty}

\begin{abstract}
We study \emph{regular spherical actions}, namely, actions of a connected compact Lie group $G$ on a manifold $M$ with equidimensional orbits, some of which are spheres.
These spherical orbits are diffeomorphic to homogeneous spaces $G/H$, and there is a list of nine possible families of pairs $(G,H)$. 

The normalizer $N$ of $H$ in $G$ plays an important role in our calculations, as does the Weyl group  $W=N/H$.
In fact, there are three possibilities for  $W$, namely: $ \Z_2, S^1$ and $S^3$.

The goal of this paper is to relate the cohomology of $M$ to that of the orbit space $M/G$. We obtain a particularly simple relationship in the first  case, $W=\mathbb{Z}_2$.

For the remaining two cases, $W=S^1$ and $W=S^3$, it is already known that the cohomology of $M$  and that of the orbit space $M/G$ are related by the Gysin sequence
$$
\cdots \to \coho H k {M/G} \to \coho H k M \to  \coho H {k-d}  {M/G}  \xrightarrow{\delta}  \coho H {k+1} {M/G}  \to \cdots.
$$
In this work, we provide a direct construction of this sequence, avoiding the use of spectral sequences. Moreover, we give an  explicit description of the connecting homomorphism as multiplication by a power of the Euler class $\epsilon \in \coho H{\ell+1}  {M/G}$ associated with the  action of  the  Weyl group $W=S^\ell$ on the fixed-point set $M^H$.

\end{abstract}

\maketitle

\vspace{-1cm}

\tableofcontents

\vspace{-1cm}

\section{Introduction}
Given a smooth $S^1$-action on a manifold $M$, there is a Gysin sequence relating the cohomology of $M$ to that of the orbit space $M/S^1$:
$$
\cdots \to \coho H i {M/S^1} \to \coho H i M \to \coho H {i-1} {M/S^1,M^{S^1}} \to \coho H {i+1} {M/S^1} \to \cdots
$$
(cf. \cite[pag. 161]{Bredon}). For an $S^3$-action, we have a more sophisticated Gysin sequence.
$$
\xymatrix@C=1.3em{
\cdots \ar[r] &
\coho{H}{i}{M}\ar[r] &
\coho{H}{i-3}{M/S^3,\Sigma/S^3}   \oplus  \left( \coho{H}{i-2}{M^{S^1}} \right)^{-{\Z}_{2}}
\ar[r] &\coho{H}{i+1}{M/S^3} \ar[r]
 &  \coho{H}{i+1}{M}\ar[r] & \cdots, &
}
$$
see \cite{RPSA} for details.

We take this one step further by aiming to extend these results to regular spherical actions. By this, we mean actions of a connected compact Lie group $G$ on a connected  manifold $M$ in which all orbits have the same dimension $d$ and at least one orbit is a sphere.
The Slice Theorem tells us that the family of spherical orbits is indeed an open dense subset of $M$.

 In future work, we will study the non-regular case, where orbits of different dimensions may occur, including, for example, actions with fixed points.

Notice that some orbits need not be spheres, for instance, real projective spaces may occur (see Example~\ref{Ejemplo}). In any case, the orbits are quotients of $S^d$ by finite subgroups.
The orbit space $M/G$ is not necessarily a manifold, but may instead be an orbifold. Consequently, the canonical projection $\pi \colon M \to M/G$ need not be a principal bundle, and hence not all the standard tools available in the principal-bundle setting can be applied.

Our goal is to establish a relationship between the cohomology $\coho H * M$ of the manifold and the cohomology $\coho H * {M/G}$ of the orbit space $M/G$. Both cohomologies are singular cohomologies with real coefficients.
The former can be computed using the de Rham complex $\coho \Om * M$, while the latter can be computed by means of the complex $\coho \Om * {M/G}$ of basic differential forms (cf. \cite{Ver}).


Given such a regular spherical action, there exists a connected closed Lie subgroup $H \subset G$ such that $G/H$ is the sphere $S^d$ and, for every point $x \in M$, the identity component of the isotropy subgroup $G_x$ is conjugate to $H$ in $G$.
Thus, there exists  an effective transitive action  of $G$ on $S^d$. It was proved in \cite{MontSam} that there are only nine distinct families of pairs $(G,H)$ satisfying this property. They are listed in \tabref{devito}.

The normalizer $N=N_G(H)$ of $H$ in $G$ plays an important role. In fact, there are three possibilities for the Weyl group
$
W=N/H,
$
namely :
$
 \mathbb{Z}_2, S^1 \text{ and } S^3.$
 
In the first case, the Gysin sequence splits, and we obtain explicit formulas for the cohomology of $M$
 (see \eqref{Z2}).

The last two cases, corresponding to $W=S^1$ and $W=S^3$, give rise to the Gysin sequence
\begin{equation}\label{GScla}
\cdots \to \coho H k {M/G} \to \coho H k M \to \coho H {k-d} {M/G,M^G} \xrightarrow{\delta} \coho H {k+1} {M/G} \to \cdots,
\end{equation}	
where the connecting morphism is given by multiplication by $\epsilon^{(d+1)/\ell}$. Here, $\epsilon \in \coho H {\ell+1} {M/G}$ denotes 
the Euler class
associated with the action of the Weyl group $W=S^\ell$ on the manifold $M^H$.

The existence of the Gysin sequence has been proved in  \cite[Theorem 4.5]{Razny2}. We provide a different proof, which does not rely on spectral sequences. Our approach also allows us to identify the connecting morphism explicitly.
As an application, we observe that the octonionic Hopf fibration does not arise from an action of this type, even though its fibers are spheres (see \remref{Hopf}).

As a further remark, we provide a procedure to describe all regular spherical actions $\Phi \colon G \times M \to M$ associated with an item in the list of \tabref{devito}. These actions are in one-to-one correspondence with almost-free smooth actions $\Psi \colon W \times E \to E$ as follows:
$
M = G/H \times_W E,
$
where the action $\Phi$ is given by the left multiplication of $G$ (cf. \propref{PT}). The manifold $E$ plays the role of $M^H$.

The main result of this work establishes how to compute the cohomology
of $M$ in terms of the cohomology groups $\coho H* {M/G}$ and
$\left(\coho H*{M^H}\right)^{-W}$ (see \teoref{nagusia}).
In the last paragraph of this work, we point out that the cohomology
of $M$ can also be computed solely in terms of the action $\Psi$,
more specifically, using the cohomology groups $\coho H *{E/W}$
and $\left(\coho H *{E}\right)^{-W}$.
In Cases II and III, the Euler class
$\epsilon \in \coho H {\dim W+1} {E/W}$ associated with the
action $\Psi$ appears to determine the connecting homomorphism
in the Gysin sequence that we construct.

The work is organized as follows. Section 2 is devoted to the study of regular spherical actions. In particular, we prove that a manifold supporting such an action is, in fact, a very particular twisted product. Section 3 is devoted to the study of the cohomology of this twisted product in terms of its factors.
The Gysin sequence is constructed and studied in Section 4.

Finally, in this work, we need to recall some properties of the Lie groups appearing in \tabref{devito}, as well as those of their Lie algebras. This is done in the Appendix, which can be read independently of the other sections.

\smallskip

In what follows, $M$ denotes a smooth connected manifold, and $G$ denotes a connected compact Lie group.

\section{Regular spherical  actions.}
\begin{quote}
{\em

We introduce regular spherical actions, that is, Lie group actions with equidimensional orbits whose principal orbits are spheres. We prove that the manifold is, in fact, a twisted product. We also recover some results previously obtained in \cite{Razny2}.

}\end{quote}

\begin{definition}
 A smooth action $\Phi \colon G \times M \to M$ is a {\em regular action} if all its orbits have the same dimension.
\end{definition}

The orbit-type stratification introduced in \cite[Chapter IV, Section 3]{Bredon} classifies the orbits of a regular acton into two types: principal orbits and exceptional orbits.
Any two principal orbits are diffeomorphic, since the isotropy subgroups of points on these orbits are conjugate. The union of all principal orbits forms an open dense subset of $M$. 
Every exceptional orbit is the base space of a covering space that is homeomorphic to a principal orbit.

The Slice Theorem (see for example \cite[Theorem  5.6]{Slice})  implies that the identity components of the isotropy subgroups of any two points in $M$ are conjugate. In other words, for any $x,y\in M$, we have $(G_x)_e\sim (G_y)_e$, although $G_x$ and $G_y$ need not be conjugate,
since they may have different numbers of connected components.

\begin{definition}
Let $\Phi\colon G\times M\to M$ be a regular action. A {\em principal isotropy subgroup} $H$ is the isotropy subgroup $G_{x_0}$ of a point $x_0$ lying on a principal orbit. The normalizer of $H$ in $G$ is denoted by $N_G(H)=N$. The quotient group $W=N/H$ is called the {\em Weyl group}.
\end{definition}

Any other choice of the base point $x_0$ yields a principal isotropy subgroup which is conjugate to $H$, while the Weyl group remains unchanged.

Note that the normalizer $N$, and hence $W$, acts on the fixed-point submanifold $M^H$. This action $\Phi_0 \colon W \times M^H \to M^H$ is effective because $G_{x_0}=H$. Moreover, if $H$ is connected, the action is almost free, and we have the following characterization:
\begin{equation}\label{Bat}
 M^H = \{ x \in M \mid (G_x)_e = H\},
\end{equation}
since $(G_x)_e \sim H$ implies $(G_x)_e=H$ for every $x\in M^H$.

\begin{definition}
An effective smooth regular action $\Phi\colon G\times M\to M$ is called a {\em regular spherical action} if one, and hence every, principal orbit is diffeomorphic to the sphere $S^d$.

\end{definition}

Notice that
 the natural projection $G/{H_e} \to G/H$ is an $H/H_e$-covering, where $H_e$ is the identity component of $H$.
The base space $G/H$  is diffeomorphic to $S^d$. When $d\ge 2$ the group $H=H_e$ is connected. 
The case $d=0$ is not possible since $G$ is connected.  For the case $d=1$ we have the following result

\begin{proposition}
Let $\Phi \colon G \times M \to M$ be a regular spherical action whose principal orbits are diffeomorphic to $S^1$. Then $G=S^1$.

\end{proposition}

\begin{proof}  Let $H$ be a  principal isotropy subgroup.
The canonical projection
$
G/H_e \longrightarrow G/H
$
is a covering map. Since $G/H \cong S^1$ and $G/H_e$ is connected and compact, we conclude that $G/H_e \cong S^1$.

We prove below that $H_e$ is a normal subgroup of $G$. Hence,
$
H_e = (G_x)_e
$
for every $x$ belonging to a principal orbit. Thus $H_e$ acts trivially on the union of principal orbits, which is a dense set in $M$. Therefore, $H_e$ acts trivially on all of $M$. Since the $G$-action is effective, it follows that
$
H_e = \{e\}.
$
This gives
$
G = S^1.
$

Let us prove the normality of $H_e$. Compactness of $G$ gives a $G$-invariant Riemannian metric $\mu$ on $G/H_e$. The Myers-Steenrod theorem
establishes that  $\Iso(G/{H_e},\mu)$ is a compact Lie group,  whose dimension is either $0$ or $1$
(see \cite[Theorem 9]{MontgomeryZippin},  \cite[Theorem 10]{MyersSteenrod}).
The homomorphism induced by the left action of $G$ on $G/H_e$,
$
\chi \colon G \to \Iso(G/{H_e},\mu),
$
cannot be constant since $G \neq {H_e}$. Hence,
$
\dim \Iso(G/{H_e},\mu)=1,
$
$\chi $ is a surjection 
and therefore
\[
\dim \Ker \chi
= \dim G - 1
= \dim G - \dim(G/{H_e})
= \dim {H_e}.
\]

Since ${H_e}$ is connected, it suffices to prove that $\Ker \chi \subset {H_e}$ to get $\Ker \chi ={H_e}$ which implies the normality of ${H_e}$.

If $g \in \ker \chi$ we have $g \cdot kH_e =kH_e$ for each $k\in G$. In particular $gH_e =g \cdot eH_e =eH_e =H_e$ which gives $g\in H_e$ and then $ \ker \chi \subset H_e$. \qedhere
\end{proof}

Since the principal orbits are circles,  the cohomology of $M$ can be computed from the cohomology of the orbit space $M/G$ using the classical Gysin sequence associated with a circle action, as recalled in the introduction.
\textbf{From now on, we assume that the spheres appearing in a regular spherical action have dimension at least $2$.}
 Hence, the principal isotropy subgroup $H$ is connected.

\bigskip
%
%

The case of a regular action yields a very explicit presentation of the manifold $M$ as a twisted product. This result is inspired by \cite[Chapter~II, Theorem~5.9]{Bredon} with some significant differences. Here, we work in the smooth category, and although all orbits have the same dimension, they are not necessarily of the same orbit type (see also \cite[Theorem 3.3]{Razny2} for another proof).

\begin{proposition}\label{PT}
Let $\Phi \colon G \times M \to M$ be a regular action of a compact connected Lie group.  Assume that the principal isotropy subgroup $H$ is connected.
Then, $M$ is $G$-equivarently diffeomorphic to the twisted product 
$$G \times_N M^H = G/H\times_W M^H,$$ where $N = N_G(H)$ and $W=N/H$ is the Weyl group.

Moreover, the orbit space $M/G$ is homeomorphic to the orbifold space $M^H/W$ of the action $\Phi_0$.
\end{proposition}
\begin{proof}
Let us consider the map
$
f \colon G \times_N M^H \to M
$
defined by
$$
f([g,z] )= \Phi(g,z) =  g \cdot z,
$$
where an element of a twisted product is denoted by $[-,-]$.
\begin{itemize}[itemsep=0.5em]

\item {\em $f$ is $G$-equivariant}: $
g' \cdot f([g,z] )= g' \cdot (g \cdot z) = g'g \cdot z = f([g'g,z]) = f (g'\cdot [g,z]).
$

\item {\em $f$ is a monomorphism.}  Assume $[g,z], [g',z']  \in G\times_H M^H$ with $f([g,z] ) = f([g',z'] )$. Condition \eqref{Bat} gives
$
H = (G_z)_e
= (G_{z'})_e$. Now, we have
\begin{align*}
f([g,z] ) = f([g',z'] ) \Rightarrow& g\cdot z = g'\cdot z' \Rightarrow g G_z g^{-1} = g'G_{z'}g'^{-1}
\Rightarrow g (G_z)_e g^{-1} = g' (G_{z'})_e g'^{-1}\\
\Rightarrow&  g H g^{-1} = g' Hg'^{-1}  \Rightarrow g'^{-1} g \in N,
\end{align*}
and therefore
$
[g,z] = [g g^{-1} g' , g'^{-1} g z]  =[g',z'].
$

\item {\em $f$ is an epimorphism.} Let $x\in M$. There exists $g\in G$ with 
$
H = g(G_x)_e g^{-1} =   (G_{g\cdot x})_e   
$
and therefore $ g\cdot x \in M^H$. So, $[g^{-1}, g\cdot x] \in G \times_N M^H $ verifies 
$
f([g^{-1}, g\cdot x] ) = x.
$
\item  {\em $f$ is  smooth } Let us consider the commutative diagram:
\begin{equation}\label{diag2}
\xymatrix{
  G\times M^H  \ar[d]^{\pi}  \ar[rd]^F  \\
  G\times_N M^H  \ar[r]^-f & M
}
\end{equation}
where $\pi$ is the canonical projection and $F(g,x) = g\cdot x$. The map $\pi$ is a $N$-principal bundle and therefore a submersion. So, the map $f$ is smooth since $F$ is smooth.

\item {\em $f$ is a submersion}. Since the map $F$ is $G$-equivariant it suffices to show that the map 
$$T_{(e,x)} F \colon T_eG \times T_xM^H \to T_xM$$ is surjective for $x\in M^H$.

By the Slice Theorem the point $x$ has a $G$-invariant neighborhood of the form $G\times_{G_x} \R^m$ where $G_x$ acts orthogonally on $\R^m$. Since the question is local, we can identify  $M$ with the twisted product $G\times_{G_x} \R^m$ and $x$ with the point $[e,0]$.

Let us compute $M^H$ in this setting. First, observe that for any $[g,z] \in G\times_{G_x} \R^m$ we have $G_[g,z] = g  (G_x)_z  g^{-1}$.
Therefore, since the action is regular, we obtain
$$
\dim H =  \dim G_{[g,z]} = \dim g (G_x)_zg^{-1} = \dim  (G_x)_z \leq \dim G_x =\dim H.
$$
Hence, equality holds throughout, and in particular,
$((G_x)_z)_e = (G_x)_e =H,
$
where the last equality follows from $x\in M^H$, which gives $(G_x)_e=H$.
Consequently, we obtain 
\begin{eqnarray*}
M^H &= &\{[g,z] \in G \times_{G_x} \R^m \mid (G_{[g,z]})_e =H\} =  \{[g,z] \in G\times_{G_x} \R^m \mid g ((G_x)_z)_e g^{-1}  =H\}\\
&=&
 \{[g,z] \in G\times_{G_x} \R^m \mid g H g^{-1}  =H\} 
= N \times_{G_x} \R^m 
\end{eqnarray*}
(cf. \eqref{Bat}).
We need to show that the map 
$$
F \colon G \times (N\times_{G_x} \R^m) \to G\times_{G_x} \R^m,
$$
given by $F(k,[g,z]) = [kg,z]$
is a submersion at the point $(e,[e,0])$,
Consider  the commutative diagram
\[
\xymatrix{
 G \times (N\times\R^m)  \ar[r]^-{\tau \times \id}  \ar[d]^{\id \times \pi}  & G \times \R^m  \ar[d]^\pi\\
 G \times (N\times_{G_x} \R^m)   \ar[r]^-{F }  & G\times_{G_x} \R^m
}
\]
where $\tau(g,k) = g\cdot k$ and $\pi(g,x) = [g,x]$. Since $\pi$ is a bundle map, and therefore a submersion, then it suffices  to verify that the top map is a submersion, or more precisely that $\tau$ is a submersion. This follows from the fact that, for each $n\in N$, the restriction $\tau \colon G \times  \{n\}  \to G$ is a diffeomorphism.

\item {\em $f$ is a diffeomorphism.}
Since $f$ is a homeomorphism, we have $\dim M = \dim(M \times_N M^H)$.
As $f$ is also a submersion, it is necessarily a local diffeomorphism,
and hence a diffeomorphism.

\item {\em Orbit space}. Since $f$ is a $G$-equivariant diffeomorphism then it induces the homeomorphism
$
\overline f \colon (M\times_N M^H)/G \to M/G.
$
It suffices to find a homeomorphism between $(M\times_N M^H)/G $ and $M^H/N$.
This is done for example in \cite[Chapter  II, Theorem  3.3]{Bredon}. The equality $M^H/N = M^H/(N/H) =M^H/W$ follows from  the fact that $H$ acts trivially on $M^H$  and that $H$ is normal in $N$. The quotient $M^H/W$ is an orbifold since the action $\Phi_0 \colon W \times M^H \to M^H$ is almost free.
\qedhere
\end{itemize}
\end{proof}

If the group $H$ is not connected, we obtain the same description of $M$ as a twisted product
$G\times_{N} M^{H_e}$
where $N=N_G(H_e)$ is the normalizer of the identity component $H_e$ in $G$. In this work, the group $H$ is connected.

\medskip

As we will see in the following lemma, the effectiveness of a regular spherical action implies that the natural action of $G$ on $G/H$ is effective. Hence, there are only finitely many possibilities for the pair $(G,H)$, as listed in \tabref{devito}.

\begin{lemma}
Let $\Phi \colon G \times M \to M$ be a regular spherical action, and let $H$ be a principal isotropy subgroup. Then the left action of $G$ on $G/H$ is effective.
\end{lemma}
\begin{proof}
Suppose that the left action of $G$ on $G/H$ is not effective. Then there exists an element $g \in G$, with $g \neq e$, such that
$
gkH = kH
$
for every $k \in G$. In particular,
$
g \in kHk^{-1}
$
for every $k \in G$.

Let $y \in M$ be a point of a principal orbit. Then
$
(G_y)_e = a H a^{-1}
$
for some $a \in G$. Hence,
$
g \in (G_y)_e,
$
and therefore $g \cdot y = y$. Since this holds for every point on a principal orbit, and the union of principal orbits is dense in $M$
\cite[Chapter  IV - Theorem  3.1]{Bredon}, it follows that $g$ fixes every point of $M$.

This contradicts the assumption that the action $\Phi$ is effective, since $g \neq e$. Therefore, the action of $G$ on $G/H$ must be effective.\end{proof}

\begin{parrafo}{\bf Classification of $(G,H)$}.
Let $G$ be a connected compact Lie group and $H$ a connected closed Lie subgroup such that $G/H$ is a standard sphere $S^d$, $d\geq 2$. We also assume that $G$ acts effectively on the sphere. The classification of such pairs $(G,H)$ was obtained by  \cite{Borel} and  \cite{MontSam}.
We find in \cite{DeVito} the following   table\footnote{The last two columns do not appear in the references cited above, but the corresponding results are scattered throughout the literature. We developed some of the proofs needed for this work in the Appendix.}

    \smallskip
   \begin{table}[H]
      \renewcommand{\arraystretch}{1.5}.
    \centering
    \begin{tabular}{|c|c|c|c|c|c|}
        \hline
        & $G$ & $H$ & $G/H$ & $N_G(H)$ & $Weyl \ group$ \\
        \hline
        
        \hline
        I-a &
        $SO(n+1)$ &
        $SO(n)$ &
        $S^n$ &
        $O(n)$ &
        $\Z_2$
        \\
        \hline
        I-b &
        $Sp(n+1)\cdot Sp(1)$ &
        $Sp(n)\cdot Sp(1)$ &
        $S^{4n+3}$ &
        $Sp(n)\cdot Sp(1)\times \Z_2$ &
        $\Z_2$
        \\
        \hline
        I-c &
        $G_2$ &
        $SU(3)$ &
        $S^6$ &
        $SU(3)\rtimes \Z_2$ &
        $\Z_2$
        \\
        \hline
        I-d &
        $Spin(7)$ &
        $G_2$ &
        $S^7$ &
        $G_2 \times \Z_2$ &
        $\Z_2$
        \\
        \hline
        I-e &
        $Spin(9)$ &
        $Spin(7)$ &
        $S^{15}$ &
        $Spin(7) \times \Z_2$ &
        $\Z_2$
        \\
                \hline
        II-a &
        $SU(n+1)$ &
        $SU(n)$ &
        $S^{2n+1}$ &
        $U(n)$ &
        $U(1)$
        \\
        \hline
        II-b &
        $U(n+1)$ &
        $U(n)$ &
        $S^{2n+1}$ &
        $U(n)\times U(1)$ &
        $U(1)$
        \\
        \hline
        II-c &
        $Sp(n+1)\cdot U(1)$ &
        $Sp(n)\cdot U(1)$ &
        $S^{4n+3}$ &
        $Sp(n)\cdot U(1) \cdot U(1)$ &
        $U(1)$
         \\
        \hline
        III &
        $Sp(n+1)$ &
        $Sp(n)$ &
        $S^{4n+3}$ &
        $Sp(n)\times Sp(1)$ &
        $Sp(1)$
        \\
        \hline
    \end{tabular}
    \medskip
    \caption{Montgomery--Samelson--Borel classification of effective
    transitive actions on spheres.}
    \label{devito}
\end{table}

The notation $Sp(m)\cdot Sp(1)$ stands for $\dfrac{Sp(m)\times Sp(1)}{\Z_2}$, where $\Z_2$ is generated by the element $(-I,-1)$. Analogously for $Sp(m)\cdot U(1)$ and $Sp(m)\cdot U(1) \cdot U(1)$.

The case $(G=S^3, H=S^1)$ does not appear in this Table because the 
action of $G$ on $G/H$ is not effective. The ineffective kernel is the center $\{\pm 1\}$ of $S^3$. 
Passing to the quotient by this subgroup, we obtain the pair $(SO(3), SO(2))$.

When $G=S^1$ or $S^3$ and $H=\{1\}$, the relationship between the cohomology of $M$ and that of $M/G$ is described by the classical Gysin sequences associated to an almost free action of $S^1$ and $S^3$, respectively, as discussed in the Introduction.
As we have already seen,  we may restrict our attention to the case $\dim G/H\geq 2$. 
Therefore, we shall assume $n\geq 2$ for $G=SO(n+1), SU(n+1)$ and $n\geq 1$ in all other cases.
\end{parrafo}

\begin{parrafo}{\bf Non-spherical orbits.}\label{Ejemplo}
Consider a regular spherical action on a manifold $M$, associated with a pair $(G,H)$ from the table above. We follow the terminology of \cite{Bredon}. The principal orbits are diffeomorphic to the sphere $S^d$. The remaining orbits, the exceptional orbits, are of the form $S^d/K$, where $K\subset O(d+1)$ is a finite group.
Let us consider an example.

 Let $(G,H)$ be a pair of groups from the table above, with $G/H = S^d$. Consider the action of $G$ on the manifold
$$
\mathbb{R}P^{d+1} \# \mathbb{R}P^{d+1}
= \frac{S^d \times [0,1]}{(x,0)\sim (-x,0),\ (x,1)\sim (-x,1)},
$$
induced by the action of $G$ on $S^d$. All orbits are diffeomorphic to $S^d$, except for the two orbits at the endpoints, which are diffeomorphic to $\mathbb{R}P^d$. Following the terminology of \cite{Bredon}, the spherical orbits are called principal orbits, while the two remaining orbits are called exceptional orbits.
\end{parrafo}

\begin{parrafo}{\bf Generic example.}\label{example}
A way to find all regular spherical actions is the following.
We choose a couple $(G,H)$ on the previous table. We consider the normalizer $N$ of $H$ and 
choose an almost free smooth action of $N/H$ on a manifold $E$. This induces 
a smooth action of $N$ on $E$, which allows us to construct the twisted product
$$
M = G \times_N E  =G/H \times_W E,
$$
where $G$ acts on the left. This action is a regular spherical action, and 
and by \propref{PT}, every action of this type
 can be obtained in this way. We necessarily have $E=M^H$.

 The points of $M$ whose orbits are not spheres (exceptional orbits) are of the form $[gH,y]$, where $y$ belongs to the singular set of $M^H$ with respect to the aforementioned action.

\end{parrafo}

\section{Cohomology of a twisted product.}\label{Tw}
\begin{quote}\em
We exhibit an explicit complex of differential forms that allows us to compute the cohomology of a twisted product $G\times_N E$ in terms of the Lie algebra of $G$ and the manifold $E$. We particularize our considerations to the cases that will arise in the paper, namely those in which the Weyl group is 
 $S^1$ or $S^3$.
\end{quote}

We consider the manifold  $M = G \times_N E$, where $G$ is a connected compact Lie group, $H \subset G$ is a connected closed Lie subgroup, and $N = N_G(H)$  denotes the normalizer of $H$ in $G$. The group $N$ acts on $E$, and the action of the subgroup $H$ on $E$ is trivial. The
induced action of the Weyl group $W= N/H$ on $E$ is almost free.
The manifold $M$ is endowed with a $G$-action induced by the left action of $G$ on itself.

In what follows, we construct a rather manageable complex that allows us to compute the cohomology of $M$ in terms of $G$ and $E$. We then particularize it to the cases of interest, namely $W=S^1$ (Case II) and $W=S^3$ (Case III).

\begin{parrafo}
\label{31} {\bf Fundamental vector fields and dual differential forms.} Let us consider a basis $\mathcal B =\{u_1, \ldots,u_p\}$  of the Lie algebra $\g$ of $G$ with:
\begin{itemize}
\item[-] $\{u_1, \ldots,u_a\}$ is a basis of the Lie algebra $\h$ of $H$, and
\item[-] $\{u_{1}, \ldots,u_{a+b}\}$ is a basis of the Lie algebra $\n$ of $N$.
\end{itemize}
The associated \emph{fundamental vector field}, relative to the right action of $G$ on itself, for a given $u\in \g$, is denoted by $Y_u$. 
These vector fields are $G$-left invariant vector fields on $G$.
 
Consider a bi-invariant Riemannian metric $\nu$ on $G$, and an $\operatorname{Ad}$-invariant scalar product $\langle -,-\rangle$ on the Lie algebra $\mathfrak g$ satisfying
$
\langle u,v\rangle = \nu(Y_u,Y_v)
$
for every $u,v \in \mathfrak g$. This is possible since $G$ is compact.

From now on, we shall assume that the basis $\mathcal B$ is orthonormal.

Given $u\in \g$ we write $\gamma_u =  i_{Y_u} \nu \in  \coho \Omega 1 G$ the associated \emph{characteristic form}.
These differential forms $\gamma_\bullet$  are $G$-left invariant.
    For simplicity, we denote
$
Y_{u_i} = Y_i,
$
and 
 $
 \gamma_{u_i} = \gamma_i,
 $
for each $i\in \{1, \ldots n\}$. 


The last ingredient we need is the family $\{Z_1, \ldots, Z_{a+b}\}$ of fundamental vector fields on $E$ defined by the action of $N$ on $E$. 
Since the group H acts trivially on E, we have that 
$
Z_1=\cdots = Z_a =0,
$
 and also that the vector fields $\{Z_{a+1}, \ldots, Z_{a+b}\}$ are, indeed, the fundamental vector fields of the almost free action of the Weyl group W on E.
 We assume that $E$ is endowed with an $W$-invariant Riemannian metric.
We denote by $\{\zeta_{a+1}, \ldots \zeta_{a+b}\}\subset \coho \Omega1 E$  the family of differential forms dual to the fundamental vector fields $\{Z_{a+1}, \ldots Z_{a+b}\}$  with respect to the invariant metric. 

In what follows, we denote by $\Omega^*_{hor}(E)$ the space of \emph{horizontal forms}, that is,
$$
\Omega^*_{hor}(E)
=
\left\{
\eta\in\Omega^*(E)\mid
i_{Z_k}\eta=0 \hbox{ for }  a+1\leq k\leq a+b
\right\}.
$$
The cohomology $\coho H *M$ can be computed using the complex of $G$-invariant differential forms on $M$:
$$
\coho \BOm *M
=
\{ \omega \in \coho \Om *M \mid L_{Y_i}\omega = 0,\ 1\leq i\leq p\},
$$
see \cite[Theorem I, Section~4.3]{GHV}. Moreover, the cohomology $\coho H*{M/G}$ of the orbit space  can also be computed using differential forms on $M$. More precisely, we compute the cohomology of  $M/G$ using the complex of basic differential forms (cf. \cite{Ver}):
$$
\coho \Om *{M/G}
=
\{ \omega \in \coho \Om *M \mid i_{Y_i}\omega = L_{Y_i}\omega = 0,\ 1\leq i\leq p\}.
$$

\begin{proposition}\label{CalTor}
Let $M=G\times_N E$ be the twisted product constructed above. Suppose that $N$ is connected. Then 
the complex $\coho \BOm *M $ is

\begin{equation}\label{Formula}
\left\{
\omega \in 
\bigwedge^*\bigl(
(\gamma_{a+1} -\zeta_{a+1}), \ldots, (\gamma_{a+b} -\zeta_{a+b}) ,\gamma_{a+b+1}, \ldots, \gamma_p
\bigr) \!\otimes \!\lau \Omega * {hor}E 
\ \middle| \
\substack{
L_{Y_j}\omega =0 \ \  \quad 1\leq j\leq a\\
L_{Y_{a+1}}\omega  = -L_{Z_{a+1} }\omega \\
\cdots\\
L_{Y_{a+b}}\omega  = -L_{Z_{a+b} }\omega 
}
\right\}
\end{equation}
\end{proposition}

\begin{proof}
Let us consider the product $G \times E$ endowed with the left $G$-action. It is straightforward to prove that
$
\coho \Om *{G \times E}
=
\bigwedge^*(\gamma_{1}, \ldots, \gamma_{p})
\otimes
\lau \Om*{hor} {G\times E},
$
and
\[
\coho \BOm *{G \times E}
=
\bigwedge^*(\gamma_1, \ldots, \gamma_p)
\otimes
\coho \Om* {E}.
\]

Here, $\lau \Om*{hor} {G\times E}$ denotes the subspace of horizontal forms with respect to the canonical projection
$
\pr \colon G\times E \to E.
$
The natural projection
$
\Pi \colon G \times E \to G\times_N E =M
$
is a $G$-equivariant bundle map. The induced map
\[
\Pi^* \colon \lau{\BOm}{*}{}{G \ib{\times}{N} E}
\to
\lau{\BOm}{*}{}{G \times E}
\]
is injective. Hence, we can identify
$
\coho{\BOm}{*}{G \ib{\times}{N} E}
$
with $\operatorname{Im}\Pi^*$, namely with the differential subcomplex of $N$-horizontal and $N$-invariant differential forms in
$
\coho \BOm * {G\times E},
$
that is,
\begin{equation*}
\coho \BOm * M =  
\left\{
\omega \in
\bigwedge^*(\gamma_{1}, \ldots, \gamma_p)
\otimes
\coho \Om * {E}
\ \middle| \
\begin{aligned}
i_{Y_j}\omega &= - i_{Z_j}\omega & 1\leq j \leq a+b,\\
L_{Y_k}\omega &= - L_{Z_k}\omega & 1\leq k \leq a+b.
\end{aligned}
\right\}.
\end{equation*}

Finally, since
$
Z_1 = \cdots = Z_a = 0,
$
the complex reduces to
\begin{equation}\label{form}
\coho \BOm * M = \left\{
\omega \in
\bigwedge^*( \gamma_{a+1}, \ldots, \gamma_p)
\otimes
\coho \Om * {E}
\ \middle| \
\begin{aligned}
i_{Y_j}\omega &= - i_{Z_j}\omega &  a+1\leq j \leq a+b,\\
L_{Y_k}\omega &= 0 & 1\leq k \leq a,\\
L_{Y_\ell}\omega &= - L_{Z_\ell}\omega &   a+1\leq \ell \leq a+b,
\end{aligned}
\right\}.
\end{equation}

Let $\omega  \in \coho \BOm * M$. We can write
$$
v = A_0 + \gamma_{a+1} \wedge A_1 +\zeta_{a+1} \ \wedge A_2 + \gamma_{a+1} \wedge \zeta_{a+1} \ \wedge A_3
$$
with
$A_\ell \in \bigwedge^*(\gamma_{a+2}, \ldots, \gamma_p) \otimes   \underbrace{\{ B \in \coho \Om* E \mid i_{Z_{a+1}}B=0\}}_{\lau \Omega* {hor, 1} E}$, for $\ell\in \{0,\ldots,3\}$.
The first condition of \eqref{form} becomes
$
A_1 +\zeta_{a+1}  \wedge A_3=  -A_2 +\gamma_{a+1} \wedge A_3
$
and therefore
$$
\omega  = A_0 +( \gamma_{a+1}   - \zeta_{a+1}  ) \wedge A_1 
 \in
 \left( \bigwedge^*( \gamma_{a+1}   - \zeta_{a+1},\gamma_{a+2} , \ldots, \gamma_p) \otimes \lau \Omega* {hor, 1} E \right) \cap 
\coho \BOm * M.
$$
 The condition $i_{Y_{a+2}} \omega = -  i_{Z_{a+2}} \omega$ does not apply  to the factor $\gamma_{a+1}   - \zeta_{a+1} $. 
Proceeding in the same way as before, we obtain that the complex $\coho \BOm*M$ is
$$
 \left( \bigwedge^*(\gamma_{a+1}   - \zeta_{a+1},\gamma_{a+2}   - \zeta_{a+2} , \gamma_{a+3}  , \ldots, \gamma_p) \otimes \lau \Omega* {hor, 1,2} E \right) \cap 
\coho \BOm * M,
$$
where 
$
\lau \Omega* {hor, 1,2} E =  \{ B \in \coho \Om* E \mid i_{Z_{a+1}}B = i_{Z_{a+2}}B =0\}.
$
Following this approach, we obtain that the complex $\coho{\BOm}{*}{M}$ is
$$
 \left( \bigwedge^*((\gamma_{a+1} -\zeta_{a+1}), \ldots, (\gamma_{a+b} -\zeta_{a+b}),\gamma_{a+b+1}  , \ldots, \gamma_p) \otimes \lau \Omega* {hor, 1,\ldots, b} E \right) \cap 
\coho \BOm * M,
$$
where 
$
\lau \Omega* {hor, 1,\ldots,b } E =  \{ B \in \coho \Om* E \mid i_{Z_{a+1}}B = \ldots =  i_{Z_{a+b}}B=0\}
=
\lau \Omega* {hor} E.
$
In other words, we get formula \eqref{Formula}.
\end{proof}
\end{parrafo}

\begin{parrafo}{\bf Particular cases.}
This formula takes a particular form in the cases of interest to us.

\smallskip
%
%
%
  
$\bullet$ Case II : $W= S^1$. Using the notations of \parref{55} 
the above calculation  \eqref{Formula} becomes 
\begin{equation}\label{FormulaBat}
\coho \BOm*{M} \,=\, 
\left\{
\begin{aligned}
\omega\, \in\,&
\bigwedge^*\bigl(
\lambda\!-\!\chi,
 \alpha_1, \ldots , \alpha_p, \beta_1, \ldots, \beta_p
\bigr)\otimes \lau \Omega * {hor}{M^H} 
\ \Biggm| \
\begin{aligned}
L_{Y_j}\omega &=0, \  1\le j\le m\\
L_{Q}\omega &=-L_{X}\omega
\end{aligned}
\end{aligned}
\right\}
\end{equation}

\smallskip

$\bullet$ Case III : $W=S^3$. Using the notations of \parref{61} 
the above calculation  \eqref{Formula} becomes
\begin{equation}\label{FormulaHiru}
\coho \BOm * M =  \left\{
\omega \! \in \!
\bigwedge^*\bigl(
\lambda_1- \chi_1,\lambda_2- \chi_2,\lambda_3-\chi_3,
\tau_{1,1} , \ldots, \tau_{4,n} \bigr)\otimes \lau \Om * {hor}{M^H} 
\! \Bigm| \!
\left.
\substack{
L_{Y_j}\omega  =   0,   \  \qquad   1 \leq j \leq m\\
L_{Q_k}\omega = \!- \!  L_{X_k}\omega,  \ 1\le k \le 3
}
\right.
\right\}
\end{equation}\smallskip

$\bullet$ When $N$ is not connected it suffices to add the condition that the differential forms are invariant by the finite group $N/N_e$. In Case I, we have $W=\Z_2$ and we get
\begin{equation}\label{FormulaBat}
\coho \BOm * M =  \left\{
\omega \in 
\bigwedge^*\bigl(
(\gamma_{a+1},  \ldots, \gamma_p
\bigr) \otimes  \lau \Omega * {}E 
\ \middle| \
L_{Y_j}\omega =0 \ 1\leq j\leq a
\right\}^W =
\left(
\coho \BOm * {G/H}  \otimes  \lau \Omega * {}E 
\right)^W 
\end{equation}
\end{parrafo}

\section{Gysin sequence}
\begin{quote}
{\em
We construct the Gysin sequence associated with a regular spherical action and study all its terms.
}
\end{quote}

We consider a regular spherical action of a connected compact Lie group $G$ on a manifold $M$, and we denote by $H$ a principal isotropy subgroup. Let $N$ be the normalizer of $H$ in $G$, and let $W$ be the associated Weyl group. The principal orbit is a sphere $S^d$ with $d\geq 2$.

The {\em Gysin term}  is the quotient complex
$$
\coho {\mathcal{G}}* M = \dfrac{ \coho \BOm * M}{\coho \Om * {M/G}}.
$$
An element of this complex is denoted by $ \llbracket \omega\rrbracket $, where $\omega \in \coho \BOm * M$.

The Gysin sequence we construct is the long exact sequence associated with the short exact sequence
$$
0 \to \coho \Om*{M/G} \to \coho \BOm * M \to \coho {\mathcal{G}}* M \to 0.
$$
In other words,
\begin{equation}\label{sucgys}
\ldots \to \coho H k {M/G} \to \coho Hk  M \to \coho H k { \coho  {\mathcal{G}}* M }  \xrightarrow{\delta}  \coho H {k+1} {M/G} \to \cdots 
\end{equation}
Our two main tasks are the computation of the cohomology
$\coho H *{\coho {\mathcal{G}} \cdot  M}$
and that of the connecting morphism $\delta$.
We focus on Cases II and III, since the other one will be treated directly (see \eqref{Z2}). We use the notation introduced in the Background sections for Cases II-III (see Appendix).

For the first task, we need the integration-along-the-fibers operator, which we now describe. 

\begin{parrafo}{\bf Integration along the fibers.}
The {\em integrating along the fibers}
$$
\Theta_M \colon \coho {\mathcal{G}} * M
\to \coho \Om {*-d}{M/G}, 
$$
which is defined as follows, depending on the case:
\begin{description}
\item[\textbf{Case II}] We have $d=2p+1$ and 
$$
\Theta_{M,3}( \llbracket \omega\rrbracket )
=
i_{W_p}  i_{V_p} \cdots i_{W_1} i_{V_1} i_{Q-X} \omega.
$$
(using the notations of \parref{55}).
\item[\textbf{Case III}] We have $d=4n+3$ and 
$$
\Theta_{M,4}( \llbracket \omega\rrbracket )
=
 i_{T_{4,n}} \cdots i_{T_{1,1}} i_{Q_3-X_3} i_{Q_2-X_2} i_{Q_1-X_1}  \omega.
$$
(using the notations of \parref{61}).
\end{description}

We now show that in both cases, the integration along the fibers is a well-defined differential operator.

\medskip

$\bullet$ Case II. We clearly have $\Theta_{M,3}( \llbracket \omega\rrbracket )=0$ when $\omega$ is a basic differential form.
The differential form 
$$ (\lambda -\chi)  \wedge  \underbrace{\alpha_1 \wedge \beta_1 \wedge \cdots \alpha_p \wedge \beta_p }_\psi \wedge \eta,
$$
with  $\eta \in \lau\Om *{hor} {M^H}$, 
belongs to $\coho \BOm * M$ if and only if 
$\eta \in \coho \Om * {M^H/N}$ (consider \eqref{FormulaBat}, (A),(B),(D)  and \eqref{Lchi}).
We then have the following vector space decomposition:
$$
\coho \BOm *M = \left\{ \omega =  (\lambda -\chi)  \wedge \psi  \wedge \eta \mid \eta \in \lau \Om *{} {M^H/N} \right\} \oplus 
\ker \Theta_{M,3}
$$
 and 
$
\Theta_{M,3}( \llbracket \omega\rrbracket )= 2  \eta \in  \coho \Om * {M^H/N =M/G}.
$
Therefore, the operator $\Theta_{M,3}$ is well defined.

 This operator  $\Theta_{M,3}$ is a differential operator if we prove that
 $$di_{W_p}  i_{V_p} \cdots i_{W_1} i_{V_1} i_{Q-X} \omega= 0 =  i_{W_p}  i_{V_p} \cdots i_{W_1} i_{V_1} i_{Q-X} d\omega$$
 for any $\omega \in  \bigwedge^*
(\lambda -\chi)  \otimes 
\bigwedge^*\bigl(
\alpha_1, \ldots \alpha_{p},\beta_1, \ldots,\beta_p
\bigr)\otimes \lau \Omega * {}{M^H/N}.
 $

The vector fields defining the integration along the fibers vanish on $M^H$, so it suffices to consider the case where 
 $
\omega 
$
is an element of $   \bigwedge^*
(\lambda -\chi)  \otimes 
\bigwedge^*\bigl(
\alpha_1, \ldots \alpha_{p},\beta_1, \ldots,\beta_p
\bigr).
$
By degree considerations, it suffices to consider the following cases, namely, those for which the codegree is less than one.
 \begin{itemize}
 \item[1)]  $\omega =\psi $.
 \item[2)] $\omega = (\lambda -\chi) \wedge \psi$.
 \item[3)] $\omega = (\lambda -\chi) \wedge i_{V_\ell} \psi$ for some\footnote{The case 
$\eta = i_{W_\ell}\psi$ for some $\ell\in \{1,\ldots,p\}$, can be treated in exactly the same way.} $\ell\in \{1,\ldots,p\}$. 
\end{itemize}
For simplicity,  we write \ $f(\omega) = i_{W_p}  i_{V_p} \cdots i_{W_1} i_{V_1} i_{Q-X} \omega$. Notice that $d\psi=0$ (cf. (D)).
Let us study the three cases.

\medskip

1)  Since $\psi$ is a cycle  of degree $2p < 2p+1$, then
$df( \psi)= 0 =f(d\psi).$

\smallskip

2) We have
$
df(\omega )= d2 =0 
$
and
$
f( d\omega )=
i_{W_p}  i_{V_p} \cdots i_{W_1} i_{V_1} \underbrace{i_{Q-X}( (d\lambda-e) \wedge\psi)}_{=_{(A), \eqref{Lchi}} 0}=0.
$

\smallskip

3) Since the degree of $\omega$ is $2p < 2p+1$, then $df(\omega) =0$.
On the other hand we have $i_{Q-X} d(\lambda-\chi) = L_{Q} \lambda  + L_X\chi =_{(A),\eqref{Lchi}} 0$ and 
therefore $$
i_{Q-X}d\omega =
  i_{Q-X} \left( d(\lambda-\chi) \wedge  i_{V_\ell} \psi 
  -
   (\lambda-\chi) \wedge L_{V_\ell} \psi \right) =_{(C)} - 2 L_{V_\ell} \psi .
$$
Consequently,
$
f(d\omega) = -2 i_{W_p}  i_{V_p} \cdots i_{W_1} i_{V_1} L_{V_\ell} \psi =0
$
since each term of $L_{V_\ell} \psi $ is not a multiple of $\psi$ (cf. (C)).
\medskip

$\bullet$ Case III.
 We clearly have $\Theta_{M,4} ( \llbracket \omega\rrbracket )=0$ when $\omega$ is a basic differential form.
The differential form 
$$ \omega =\underbrace{(\lambda_1 -\chi_1)  \wedge  (\lambda_2 -\chi_2)  \wedge  (\lambda_3 -\chi_3)}_{\theta_1} \wedge\underbrace{ \tau_{1,1} \wedge \tau_{2,1} \wedge \tau_{3,1} \wedge \tau_{4,1}  \wedge \cdots \wedge  \tau_{1,n} \wedge \tau_{2,n} \wedge \tau_{3,n} \wedge \tau_{4,n}}_{\kappa_{vol}}  \wedge \eta,
$$
with  $\eta \in \lau\Om *{hor} {M^H}$, 
belongs to $\coho \BOm * M$ if and only if 
$\eta \in \coho \Om * {M^H/N}$ (consider \eqref{FormulaHiru}, \lemref{kappae}, \parref{52} and \eqref{Lchi}).
We then have the following vector space decomposition:
$$
\coho \BOm *M = \left\{ \omega =  \theta_1  \wedge \kappa_{vol}  \wedge \eta \mid \eta \in \lau \Om *{} {M^H/N} \right\} \oplus  \ker \Theta_{M,4} 
$$
  and 
$
\Theta_{M,4}( \llbracket \omega\rrbracket )=  8 \eta \in  \coho \Om * {M^H/N} =  \coho \Om * {M/G}.
$
Therefore, the operator $\Theta_{M,4}$ is well defined.

 This operator is a differential operator if we prove 
  $$d i_{T_{4,n}} \cdots i_{T_{1,1}} i_{Q_3-X_3} i_{Q_2-X_2} i_{Q_1-X_1} \omega= 0 =  i_{T_{4,n}} \cdots i_{T_{1,1}} i_{Q_3-X_3} i_{Q_2-X_2} i_{Q_1-X_1}  d\omega$$
  for any $\omega \in  \bigwedge^*
\bigl((\lambda_1-\chi_1),(\lambda_2-\chi_2),(\lambda_3-\chi_3) \bigr)  \otimes 
\bigwedge^*\bigl(
\tau_{1,1} , \ldots, \tau_{4,n}
\bigr)\otimes \lau \Omega * {}{M^H/N}.
 $

The vector fields defining the integration along the fibers vanish on $M^H$, so it suffices to consider the case where 
 $
\omega 
$
is an element of
 $
\omega \in  \bigwedge^*
\bigl((\lambda_1-\chi_1),(\lambda_2-\chi_2),(\lambda_3-\chi_3) \bigr)  \otimes 
\bigwedge^*\bigl(
\tau_{1,1} , \ldots, \tau_{4,n}
\bigr).
$
 Notice that $d\kappa_{vol}=C_n d\kappa^n =0$ (cf. \parref{52}).
By degree considerations, it suffices to consider the following cases, namely, those for which the codegree is less than one.
 \begin{itemize}
 \item[1)]  $\omega = \theta_1 \wedge \kappa_{vol}$.
 \item[2)] $\omega = i_{Q_j -  X_j} \theta_1 \wedge \kappa_{vol}  $, with $j \in \{1,2,3\}$,
 \item[3)] $\omega = \theta_1 \wedge i_{T_{1,1} } \kappa_{vol}$\footnote{The other cases $T_{2,1}, \ldots, T_{4,n}$ 
 can be treated in exactly the same way.}
 \end{itemize}
To simplify the notation,  we write \ $f(\omega) = i_{T_{4,n}} \cdots i_{T_{1,1}} i_{Q_3-X_3} i_{Q_2-X_2} i_{Q_1-X_1} \omega$
and $\beta_k = \lambda_k-\chi_k$ for $k\in\{1,2,3\}$.
Let us study the three cases.

1) We have
$
d f( \omega )= d8 =0 .
$
On the other hand, 
\begin{align}\label{QQQd}
i_{Q_3 -X_3} i_{Q_2 - X_2}  i_{Q_1-X_1}  d\theta_1
& =
 i_{Q_1 -X_1} d\beta_1 
 +
 i_{Q_2 -X_2} d\beta_2 
 +
 i_{Q_3-X_3} d\beta_3  
 \\ \nonumber
 &=
L_{Q_1 -X_1} \beta_1 
 +
 L_{Q_2 -X_2} \beta_2 
 +
 L_{Q_3-X_3} \beta_3=_{ \eqref{22BC},(b)}0,
\end{align}
and therefore,
\begin{align*}
f(d\omega ) &=
 f(d\theta_1  \wedge\kappa_{vol})=
i_{Q_3-X_3} i_{Q_2-X_2} i_{Q_1-X_2} d\theta_1  \wedge  i_{T_{4,n}} \cdots i_{T_{1,1}}  \kappa_{vol} =0\cdot 1 =0.
\end{align*}

\smallskip

2)  Since the degree of $\omega$ is $4n+2 < 4n+3$ then
$df( \omega) = 0 $. We may assume $j=1$, as the other cases are analogous.
Since
\begin{align*}
i_{Q_3-X_3 } i_{Q_2-X_2}  i_{Q_1-X_1}d i_{Q_1-X_1 } \theta_1
&=
i_{Q_3-X_3 } i_{Q_2-X_2}   L_{Q_1-X_1} i_{Q_1-X_1 } \theta_1 
=
i_{Q_3-X_3 } i_{Q_2-X_2}   L_{Q_1-X_1}(\beta_2 \wedge \beta_3)\\
&=_{(b),\eqref{22BC}}
i_{Q_3-X_3 } i_{Q_2-X_2}  (-(\lambda_3+\chi_3) \wedge \beta_3 + \beta_2 \wedge(\lambda_2+\chi_2)  ) \\
&=
2 i_{Q_3-X_3 } i_{Q_2-X_2}  (-(\lambda_3 \wedge \chi_3)  + (\lambda_2\wedge \chi_2)  )=0,
\end{align*}
then
$
f( d\omega) =
f (d i_{Q_1- X_1} \theta_1   \wedge\kappa_{vol})
=
i_{Q_3-X_3 } i_{Q_2-X_2}  i_{Q_1-X_1}d i_{Q_1-X_1 } \theta_1   \wedge\ i_{T_{4,n}} \cdots i_{T_{1,1}} \kappa_{vol}  =0\cdot 1 =0.
$

\smallskip

3) Since the degree of $\omega$ is $4n+2 < 4n+3$, then $d f(\omega) =0$. Since
\begin{align*}
 i_{T_{4,n}} \cdots i_{T_{1,1}} d i_{T_{1,1}}  \kappa_{vol} 
 &=
  i_{T_{4,n}} \cdots i_{T_{1,2}}  L_{T_{1,1}}  i_{T_{1,1}}   \kappa_{vol} 
  =
   - i_{T_{4,n}} \cdots   i_{T_{1,3}}  L_{T_{1,1}} i_{T_{1,2}}  i_{T_{1,1}}   \kappa_{vol} 
    \\
    &+
      i_{T_{4,n}} \cdots i_{[T_{1,2},T_{1,1}]}  i_{T_{1,1}}   \kappa_{vol} 
      =_{(vi)}    - i_{T_{4,n}} \cdots   i_{T_{1,3}}  L_{T_{1,1}} i_{T_{1,2}}  i_{T_{1,1}}   \kappa_{vol} 
      = \cdots \\
      &= - L_{T_{1,1}} i_{T_{4,n}} \cdots   i_{T_{1,3}}   i_{T_{1,2}}  i_{T_{1,1}}   \kappa_{vol} = L_{T_{1,1}} 1=0,
\end{align*}
then
\begin{align*}
f (d \omega) &= 
f(d (\theta_1 \wedge i_{T_{1,1}} \kappa_{vol} ) )=
i_{T_{4,n}} \cdots i_{T_{1,1}} (  (\underbrace{i_{Q_3 -X_3} i_{Q_2 - X_2}  i_{Q_1-X_2}  d\theta_1}_{=_\eqref{QQQd} 0} ) \wedge 
i_{T_{1,1}}\kappa_{vol} )  \\
  & - i_{Q_3 -X_3} i_{Q_2 - X_2}  i_{Q_1-X_2}  \theta_1\wedge \underbrace{i_{T_{4,n}} \cdots i_{T_{1,1}} d i_{T_{1,1}}  \kappa_{vol} }_{=0} =0. 
\end{align*}

\smallskip

We have just proved that the operator  $\Theta_M$ is a differential operator for both cases.
It also is  natural with respect to invariant open subsets. More precisely, if $U\subset M$ is an invariant open subset, then the following diagram commutes:
\begin{equation*}\label{Theta}
\xymatrix{
\coho {\mathcal G} * M \ar[r] \ar[d]^{\Theta_M}& \coho {\mathcal G} * U \ar[d]^{\Theta_U}\\
\coho \Om {*-b} {M/G} \ar[r] & \coho \Om {*-b}{U/G}
}
\end{equation*}
where horizontal maps are restrictions.
\end{parrafo}

The next step is to develop the classical Mayer-Vietoris argument. We begin with the following elementary step.

\begin{lemma}\label{MV1}  Let $\{U,V\}$ be a $G$-invariant open covering  of $M$. If  the operators 
$
\Theta_{U\cap V}, \Theta_U, \Theta_V $  are quasi-isomorphisms, then the operator $ \Theta_M $ is also a quasi-isomorphism.
\end{lemma}
\begin{proof}
We have the following commutative diagram, in which the horizontal maps are restriction maps.

$$
\xymatrix{
0 \ar[r]  &\coho \Om * {\mathcal G  (M) } \ar[r]  \ar[d]^{\Theta_M} & 
\coho \Om * {\mathcal G ( U) } \oplus  \coho \Om * {\mathcal G ( V) }
\ar[r] \ar[d]^{\Theta_U \oplus \Theta_V} &  \coho \Om * {\mathcal G ( U\cap V) } \ar[r] \ar[d]^{\Theta_{U\cap V} } & 0 \\
0\ar[r] &\coho \Om * {M/G} \ar[r] &  \coho \Om * { U/G}  \oplus  \coho \Om * { V/G} \ar[r]
&  \coho \Om * {(U\cap V)/G} \ar[r] &0.
}
$$
To prove that the short horizontal sequences are exact, it suffices to show that the two horizontal arrows on the right are surjective.
Since the Lie group is compact, there exists a $G$-invariant partition of unity $\{f,g\}$ subordinate to the covering $\{U,V\}$. So we have
\begin{align*}
\omega \in \coho \Om * {U\cap V} \Longrightarrow &( f \omega,g\omega)  \in
 \coho \Om * {U} \oplus  \coho \Om * {V} \\
\omega \in \coho \Om * {(U\cap V)/G} \Longrightarrow& ( f \omega,g\omega)  \in
 \coho \Om * {U/G} \oplus  \coho \Om * {V/G} 
 \end{align*}
since   $f,g \in \coho \Om 0 {M/G}$.
The result follows from the Five Lemma.
\end{proof}

Next, we study the case where the manifold $M$ is an orbit $G(y)$. Notice that $G(y)/G$ is an atom.

\begin{lemma}\label{MVOrb}
For each orbit $G(y) \subset M$ the operator $\Theta_{G(y)} \colon \coho {\mathcal{G}} * {G(y)}
\to \coho \Om {*-d}{G(y)/G },$ is a quasi-isomorphism.
\end{lemma}
\begin{proof}
Notice that we are in Case II or III where the normalizer $N$  is connected.
Since $M = G\times_N M^H$ we may assume, without loss of generality, that   $y\in M^H$.
Notice that  $G_y = N_y$.

The natural projection $\underbrace{G/H}_{S^d} \longrightarrow \underbrace{G/G_y}_{G(y)}$ is a $G_y/H$-covering. 
Since  $G_y/H$ is included in the connected group $N/H$, its action on $G/H$, as the restriction of the $N/H$-action,  preserves the orientation.
Therefore
$\coho H *{G(y)} =\coho H *{S^d}^{G_y/H}= \coho H *{S^d}$.
 Notice that this is not always the case in Case I.

Consider the diagram
$$
\xymatrix{
0 \ar[r] &  \coho  \Omega * {G(y)/G =pt} \ar[r] & \coho \BOm * {G(y)}   \ar[r]^\pi  &\ar[r] \ar[d]^-{\Theta_{G(y)} }\coho {\mathcal{G}} *{G(y)} & 0.\\
&&& \coho \Om {*-d} {G(y)/G =pt}.&
}
$$
It suffices to prove that $\Theta_{G(y)}^* \colon \coho H {d}  { \mathcal{G} (G(y))} =\R  \to \coho H 0 {pt}=\R$ is an isomorphism. Thus, we need to find a cycle $\llbracket \omega \rrbracket  \in  \coho { \mathcal{G}}  {d}  {G(y)} $ such that $\Theta_{G(y)} (\llbracket \omega \rrbracket )\ne 0$.
The differential form $\omega = \theta\wedge  \kappa_{vol}$ 
belongs to $\coho \BOm {d=4n+3} M$ by \lemref{kappae} and \parref{52}. Moreover,
 $\Theta_M(\llbracket \omega\rrbracket) =8$ and, 
 $$d\omega =_{\ref{kappae}, \ \ref{52}} ( \kappa -e) \wedge  \kappa_{vol}
 =
  \ \kappa \wedge \kappa_{vol} -   \ e \wedge \kappa_{vol} = - e \wedge  \kappa_{vol},
 $$
 since $ \kappa \wedge \kappa_{vol} \in \bigwedge^{4n+4}  \left(  \tau_{1,1}, \ldots , \tau_{4,n}\right) = \{0\}$. 
  The Euler form $e$ belongs to $\coho \Om 4 {N(y)/N=\{pt\}} =0$ and therefore $d\omega=0$.
\end{proof}

\begin{proposition}\label{thetaqi}
Let $G$ be a connected compact Lie group with a regular spherical action on a manifold $M$ in Cases II or III. Then,
$\Theta_M \hbox{ is a quasi-isomorphism.}
$
\end{proposition}
\begin{proof}

We associate a numerical invariant to every $G$-invariant open subset $U$ of $M$, namely its depth, denoted by $\depth U$. The key point is that we will be able to find a covering of $M$ by $G$-equivariant open subsets whose elements, up to a $G$-invariant retraction,  have strictly smaller depth than $M$. An inductive argument, together with the previous lemmas, will then allow us to obtain the desired result.

The depth $d(U)$ is  determined by the number of isotropy subgroups, up to conjugacy, contained in $U$. To give a precise definition, we use the notion of orbit type introduced in  \cite[Chapter IV, Section 3]{Bredon} (see also \cite[Section  1.2]{AndreaFranco}).

The equivalence relation $x\sim y$ defined by the condition that  ``$G_x \hbox{ congruent to } G_y$", induces a stratification $\mathcal S_M$ on $M$. Its elements are $G$-invariant submanifolds, called {\em strata}. This family is endowed with the partial order
$
S_1 \preceq S_2 \Longleftrightarrow  S_1 \subset \overline{S_2}.
$
The {\em depth} of the stratification is defined by:
$$
\depth{M} = \max \{ p \in \N \mid \hbox{ there exists a chain } S_0 \prec S_1 \prec \ldots \prec S_p\}.
$$
The Slice Theorem implies that this number is finite.	
In particular, if $\depth M=0$, then $\mathcal S_M$ consists of a single stratum, namely $M$. In this case, the canonical projection
$
\pi\colon M\to M/G
$
is a $G/H$-bundle over the manifold $M/G$.

Notice that the closed strata of $\mathcal S_M$ are precisely  the minimal strata. Each on them admits  a $G$-invariant tubular neighborhood $\tau_S \colon T_S \to S$. Then $T_S\menos S$ is a $G$-invariant open subset of $M$ satisfying $\depth
{T_S\menos S} < \depth {T_S}  \leq \depth M$. Moreover, the action of $G$ in $S$ is regular, spherical and satisfies $\depth S=0$.

Let $\Sigma$ be the disjoint union of closed strata. For each closed stratum $S$, fix a tubular neighborhood $T_S$ as above. These neighborhoods can be chosen to be pairwise disjoint. Let $T$ be their union. Then the $G$-invariant covering
$
\{M\setminus\Sigma,T\}
$
satisfies the properties we were looking for
\begin{itemize}
\item[i)] $\depth {M\setminus\Sigma}<\depth M$,
\item[ii)] $\depth {T\setminus\Sigma}<\depth M$, and
\item[iii)] $T$ $G$-equivariantly retracts onto $\Sigma$, and $\depth \Sigma=0$.
\end{itemize}

We now prove the claim of the Proposition: the operator  $\Theta_{M} $ is a quasi-isomorphism.
 In fact,  we prove by induction on the depth, the following assertion
\begin{center}
$\Pi(n)$ = ``$\Theta_M$ is a quasi-isomorphism for each manifold $M$ \\ endowed with a $G$ regular spherical action with $\depth M\leq d$,''
\end{center}
for every $n\in \N$.

\medskip

$\bullet \ \Pi(0)$. We know that the canonical projection $\pi \colon M \to M/G$ is a $G/H$-bundle over the manifold $M/G$. For each open subset $U \subset M/G$, consider the assertion
$$
P(U) = ``\Theta_{\pi^{-1}(U)} \hbox{ is a quasi-isomorphism".}
$$

If we can apply \cite[Chapter V, Lemma 9.5]{Bredon1993}, then we obtain P(M), that is, $\Pi(0)$.  That Lemma is a nice Mayer-Vietoris argument to get global properties by checking them locally, with three simple conditions. Condition (2) is a Mayer-Vietoris type property, which is granted by  \lemref{MV1}, while condition (3), which is compatibility with disjoint unions, is straightforward. It remains to prove condition (1), which is just the local checking, that is, to check whether $P(U)$ holds for $U$.  In our case, the open subset U is contractible, and hence $\pi \colon \pi^{-1} (U) \to U$ becomes the canonical projection $\operatorname{pr} \colon U \times G/ H\to U$. Since $U$ can be contracted to a point, $P(U)$ reduces to the statement that ``$\Theta_{G/H}$ is a quasi-isomorphism''. Since $G/H$ is an orbit of the action,  this follows from \lemref{MVOrb}.

\medskip

$\bullet \ \Pi(<q) \Rightarrow \Pi(q)$. Let $M$ be a manifold endowed with a regular spherical $G$ -action  such that $\depth M = q$.
Consider the $G$-invariant open covering $\{M\setminus\Sigma,T\}$ of $M$. By properties i) and ii), we obtain the claim by applying \lemref{MV1}, provided that we prove that $\Theta_T$ is a quasi-isomorphism. By iii), it suffices to prove that $\Theta_\Sigma$ is a quasi-isomorphism. This follows from $\depth\Sigma=0$ and the previous case.
\end{proof}

This Proposition gives the Gysin sequence 
$$
\cdots \to \coho H k {M/G} \to \coho H k M \to  \coho H {k-d}  {M/G}  \xrightarrow{\delta}  \coho H {k+1} {M/G}  \to \cdots.
$$
that we can find in \cite[Theorem 4.5]{Razny2} proved in a different way. 
We arrive to the main result of this work.
Recall that if the group $\mathbb{Z}_2$ acts on a differential complex $A$ and $J$ is a generator of $\mathbb{Z}_2$, then we write
$$
A^{-\mathbb{Z}_2} = \{a \in A \mid J(a) = -a\}.
$$

\begin{theorem}\label{nagusia}
Let $G$ be a connected compact Lie group acting regularly and spherically on a connected manifold $M$. Let $H$ be a principal isotropy subgroup and $N$ its normalizer. The dimension of the principal orbits is $d$.
According to the three possible types of Weyl group $W=N/H$, we obtain the following relationships between the cohomology groups of $M$ and $M/G$.

\medskip

Case I ($W=\Z_2$) -  The cohomology of the Gysin term is $\coho H {*-d}{M/G)}$ when $d$ is odd, and $\left(\coho H {*-d}{M^H}\right)^W$ otherwise. The connecting morphism $\delta$ of the Gysin sequence \eqref{sucgys} vanishes, which yields the formula
\begin{equation}\label{Z2}
\coho H *M = \coho H * {M/G} \oplus
\begin{cases}
\coho H {*-d}{M/G}, & \text{if } d \text{ is odd},\\[4pt]
\bigl(\coho H{*-d}{M^H}\bigr)^{-\mathbb{Z}_2}, & \text{if } d \text{ is even}.
\end{cases}
\end{equation}

\smallskip

Case II ($W=S^1$) - Here $d=2p+1$. The cohomology of the Gysin term is $ \coho H {*-2p-1}{M/G} $. The Gysin sequence \eqref{sucgys}
becomes
\begin{equation}\label{Gys2}
\ldots \to \coho H k {M/G} \to \coho Hk  M \to  \coho H {k-2p-1}{M/G}\xrightarrow{\delta}  \coho H {k+1} {M/G} \to \cdots 
\end{equation}
where the connecting morphism $\delta$ 
is  given by multiplication by $\epsilon^{p+1}$, where $\epsilon \in \coho H 2{M/G}$ denotes the Euler class of the action action $\Phi_0 \colon W \times M^H \to M^H$.
\smallskip

Case III ($W=S^3$) - Here $d=4n+3$. The cohomology of  the Gysin term is $ \coho H {*-4n-3}{M/G} $. The Gysin sequence \eqref{sucgys}
becomes
\begin{equation}\label{Gys3}
\ldots \to \coho H k {M/G} \to \coho Hk  M \to  \coho H {k-4n-3}{M/G}\xrightarrow{\delta}  \coho H {k+1} {M/G} \to \cdots 
\end{equation}
where the connecting morphism $\delta$ 
is given by multiplication by $\epsilon^{n+1}$, where $\epsilon \in \coho H 4{M/G}$ denotes the Euler class of the action action $\Phi_0 \colon W \times M^H \to M^H$.

\end{theorem}
\begin{proof} We consider our three cases.

\smallskip

\noindent {\bfseries Case I.} Since $W = \Z_2$, it follows from \propref{PT}  that
$$
M = G \times_N M^H = G/H \times_{W}  M^H.
$$
Notice that $M^H/W = M/G$. 
We have seen that the generator  $J$ of $W= \Z_2$ acts on the sphere $G/H = S^d$ via the antipodal map (cf. \propref{Z2}).

\negro

The Gysin sequence arises from the short exact sequence
$$
0 \to \coho \Om * {M/G} = \coho \Om * {M^H} \to \coho \BOm * M =_{\eqref{FormulaBat}} \left(\coho \BOm *{G/H} \otimes \coho \Om * {M^H}\right)^W
\to \coho {\mathcal G} *M \to 0.
$$
Since $W$ is finite then  the canonical projection $\pi \colon G/H \times M^H \to G/H \times_{W}  M^H$ is a covering map and therefore 
\begin{align*}
\coho H * M = &
\coho H* {\coho \Om \cdot {G/H \times_{W}  M^H} }  = \coho H * {  {\coho \Om \cdot {G/H \times  M^H}}^{W}}\\
&=
 \left(\coho H * {  {\coho \Om \cdot {G/H \times  M^H}}}\right)^{W}
 =
 \left( \coho H * { G/H \times  M^H}\right)^{W}\\
 &=
  \left( \coho H * { G/H}  \otimes  \coho H * {M^H} \right)^{W}= 
    \left(  \coho H * {M^H} \right)^{W}   \oplus  \left(   {\coho H d { S^d}  \otimes \coho H {*-d}  {M^H}}\right)^{W}
    \\
    =&
  \coho H * {M^H/W}   \oplus  \left(   {\coho H d { S^d}  \otimes \coho H {*-d}  {M^H}}\right)^{W}
  =
    \coho H * {M/G}   \oplus  \left(   {\coho H d { S^d}  \otimes \coho H {*-d}  {M^H}}\right)^{W}
\end{align*}
The Gysin sequence splits, that is $\delta=0$, and the cohomology of the Gysin term is given by,
$$
\coho H * {\lau {\mathcal G} \cdot V M}
=
\left\{
\begin{array}{ll}
 \left(\coho H {*-d}  {M^H}\right)^{W} = \left(\coho H {*-d}  {M^H/W}\right) = \coho H {*-d}  {M/G} & \hbox{if  $d$ odd},\\
\left(\coho H {*-d} {M^H}\right)^{-W} = \left(\coho H {*-d} {M^H}\right)^{-W} & \hbox{if $d$ even},
\end{array}
\right.
$$
where we use the notation
$
\left(\coho H * {M^H}\right)^{-\Z_2} = \{ \omega \in \coho H * {M^H} \mid J^* \omega = -\omega\}$. This gives \eqref{Z2}.
In Cases I-a (odd dimension), I-b, I-d, and I-e, we obtain
$$
 \coho H * M =    \coho H * {M/G} \otimes \coho H *{G/H} =   \coho H * {M/G} \otimes \coho H * {S^d}.
 $$\negro
These results coincide with   \cite[Corollary 4.9]{Razny2}.
The case $G=S^3$ and $H=S^1$ was studied in \cite[Remark 5.1]{RS1}. 
\negro

\smallskip

\noindent {\bfseries Case II.} Here $S^d = S^{2p+1}$. We know from \propref{thetaqi} that the  cohomology of the Gysin term is $ \coho H {*-2p-1}{M/G} $.  It remains to compute the connecting morphism $\delta$.

Let us consider a cycle $\tau \in \coho \Om * {M/G}$. We construct a cycle $\llbracket \omega \rrbracket \in \coho { \mathcal{G}} *M$ with $\Theta_{M,3}(\llbracket \omega \rrbracket) = \tau$. The connecting morphism is given by
$
\delta ([\tau]) = d \omega.
$
The differential form 
$$
\omega = \frac 1 {2C}  \sum_{k=0}^p (\lambda-\chi) \wedge (d\lambda)^k \wedge e^{p-k} \wedge \tau
$$
belongs to $ \coho \BOm *M $ (cf. \eqref{FormulaBat}, (A), (D)  and \eqref{Lchi}). 
Recall that $e\in \coho \Om 2 {M/G}$ denotes  the Euler form of the action of the Weyl group $W$ on $M^H$ (cf. \parref{Euler1}).
By degree considerations and (D) we have
$$
\Theta_{M,3}(\llbracket\omega \rrbracket)  =  \frac 1 {2C}\sum_{k=0}^p \Theta_{M,3}(\llbracket  (\lambda-\chi) \wedge (d\lambda)^k \wedge e^{p-k} \wedge \tau\rrbracket ) =
 \frac 1 {2C} \Theta_{M,3}( \llbracket (\lambda-\chi) \wedge (d\lambda)^p  \wedge \tau\rrbracket )
=  \tau.
$$
Finally
\begin{align*}
d\omega &=   \frac 1 {2C} \sum_{k=0}^p d((\lambda-\chi) \wedge (d\lambda)^k \wedge e^{p-k} \wedge \tau) \\
&=
 \frac 1 {2C} \sum_{k=0}^p (d\lambda)^{k+1} \wedge e^{p-k} \wedge \tau
- \frac 1 {2C}  \sum_{k=0}^p (d\lambda)^k \wedge e^{p-k+1} \wedge \tau
=
-  \frac 1{2C} e^{p+1} \wedge \tau \in \coho \Om *{M/G},
\end{align*}
which gives
$$
\delta([\tau]) = -  \frac 1 {2C}  \ \epsilon ^{p+1} \wedge [\tau]. 
$$
Recall that $\epsilon \in \coho H 2 {M/G}$ is the Euler class of the Weyl-group action on the fixed-point set $M^H$.

In order to  obtain the announced connecting morphism, it suffices to replace the integration $\Theta_{M,3}$
 by its composition with multiplication by $-\tfrac1{2C}$.

\smallskip

\noindent {\bfseries Case III.} Here $S^d = S^{4n+3}$. We know from \propref{thetaqi} that the  cohomology of the Gysin term is $ \coho H {*-4n-3}{M/G} $. It remains to compute the connecting morphism $\delta$.

Let us consider a cycle $\tau \in \coho \Om *{M/G} = \coho \Om * {M^H/N}$. We construct a cycle $\llbracket \omega \rrbracket \in \coho { \mathcal{G}}{* }M$ with $\Theta_M(\llbracket \omega \rrbracket) = \tau$. The connecting morphism is given by
$
\delta ([\tau]) = d \omega.
$

 From \lemref{kappae}, there exists $\theta \in \coho \BOm 3 M$ with $\kappa = e + d\theta$,
and therefore
$$
\kappa_{vol} = {C_n} \cdot \kappa^{n}  = C_n \cdot e^{n} +  d \sigma,
$$
for some $\sigma \in \coho \BOm {4n-1} M$ (cf. \eqref{123}).
Recall that $e\in \coho \Om 4 {M/G}$ denotes the Euler form of the action of the Weyl group $W$ on $M^H$ (cf. \parref{kappatau}).

The differential form we are looking for is 
$$
\omega =\frac 1 8 (  \theta \wedge \kappa_{vol} \wedge \tau  + e \wedge \sigma\wedge \tau).
$$
It belongs to $ \coho \BOm {*} M $ since $\tau, e \in \coho \Om * {M/G}$,  by \lemref{kappae}  and \parref{52}. Moreover, it satisfies
\begin{align*}
\Theta_{M,4}(\llbracket \omega \rrbracket )  &=   \tfrac 1{8} \cdot  \Theta_{M,4}(\llbracket \theta\wedge \kappa_{vol}\wedge \tau \rrbracket ) +  \tfrac 1 {8} \cdot \Theta_M(\llbracket e \wedge \sigma \wedge \tau  \rrbracket ) = \tau -0 = \tau, \\
&\hbox{ since $e \wedge\tau \in \coho \Om * {M^H/N}$ and the degree of $\sigma$ is $4n-1 < 4n+3$, and}
\\
d\omega &=  \tfrac 1 {8} \left(d \theta\wedge \kappa_{vol}  \wedge \tau  + e \wedge d\sigma \wedge \tau\right)
 = \tfrac 1 {8} \left( \kappa \wedge \kappa_{vol}  \wedge \tau - e \wedge \kappa_{vol}   \wedge \tau  + e \wedge d\sigma \wedge \tau \right)\\
 &=\tfrac 1 {8} \left( -e \wedge \kappa_{vol}   \wedge \tau + e \wedge d\sigma \wedge \tau\right)
 =
- \tfrac 1 {8} \left(  e \wedge C_n \cdot e^n  \wedge \tau\right)
= 
-\tfrac{C_n}8 \cdot e^{n+1}  \wedge \tau.
\end{align*}
We have used the equality $d\theta = \kappa - e$ from \lemref{kappae} and the fact that $\kappa \wedge \kappa_{vol} \in \bigwedge^{4n+4}  \left(  \tau_{1,1}, \ldots , \tau_{4,n}\right) = \{0\}$. 
This gives
$$
\delta([\tau]) = -\tfrac{C_n}8 \cdot   [  e^{n+1}    \wedge \tau] =  -\tfrac{C_n}8 \cdot  \epsilon ^{n+1} \wedge [\tau] \in \coho H {*}{M/G}.
$$
Recall that $\epsilon \in \coho H 4 {M/G}$ is the Euler class of the Weyl-group action on the fixed-point set $M^H$.

In order to  obtain the announced connecting morphism, it suffices to replace the integration $\Theta_{M,4}$  by its composition with multiplication by $-\tfrac  {C_n}8$.
\end{proof}

\begin{remark}\label{Hopf}

\phantom .

$\bullet$
The generalized Hopf bundle $S^7 \to  S^{15} \to S^8$ does not arise from a regular spherical action.  One indication of this is that the classical Gysin sequence associated with this spherical bundle does not degenerate. 
In Case I of this work, the Gysin sequence always degenerates.
Since $M/G=S^8$, the same conclusion holds in Cases II and III. Indeed,  the Euler class $\epsilon \in \coho H 2 {S^8}$ in Case II (respectively, $\epsilon \in \coho H4{S^8}$ in Case III) vanishes, and hence $\delta=0$.

\smallskip

$\bullet$ 
The exotic term $\left(\coho H *{M^H}\right)^{-W}$ can be expressed alternatively as
$
\frac{\coho H *{M^H}}{\coho H *{M^H/W}}
=
\frac{\coho H *{M^H}}{\coho H *{M/G}}.
$

\smallskip

$\bullet$ Note that the Euler class appearing in the Gysin sequence of an oriented sphere bundle vanishes when the dimension of the sphere is even. The same phenomenon occurs for a regular spherical action.
\end{remark}

\section{Appendix}

 We collect the relevant properties of the Lie groups and their associated Lie algebras arising from a regular spherical action. Case III is treated first, as it is needed before considering Cases I-b and II-c.\negro

This Appendix is independent of the results proved in the preceding sections.

\bigskip

\centerline {
\fbox{\begin{tabular}{c}
Background on Case III\\
$(G,H) = (Sp(n+1), Sp(n))$ with $n\geq 1$.
\end{tabular}}

}

\bigskip

Recall that $G= \{ A \in \mathcal M_{n+1}(\H) \mid A^\dagger A=I\}$. 
The inclusion $H \hookrightarrow G$ is defined by $A \mapsto
A_1 = 
\left(\begin{array}{c|c}
A & 0 \\
\hline
0 & 1
\end{array}\right)$. The normalizer $N$ of $H$ in $G$ is
\begin{equation}\label{N1}
N = \left\{ \left(\begin{array}{c|c}
A & 0 \\
\hline
0 & q
\end{array}\right)
\mid A \in H, q \in Sp(1) \right\} \cong H\times Sp(1).
\end{equation}
Let us prove this. Consider the natural action of $Sp(n+1)$ on $\H^{n+1} = \H^n \oplus \H$ given by $(B,v)\longmapsto Bv$. 
Since $
(\H^n\oplus \H)^H = \{(0,...,0,h) \mid  h\in \H\} = \H_1\cong \H
$
the normalizer $N$ preserves the subspace $\H_1$. 
Indeed, if $B\in N$ and $v\in\H_1$ then, for every $A\in H$, $B^{-1}A_1 B=A'_1$ for some $A'\in H$, and hence
$A_1 Bv = BA'_1v = Bv$. Thus, $B(\H) \subset  \left(  \H^n \oplus \H\right)^H  = \H_1$ as claimed. 

Since $B$ is unitary, it also preserves the orthogonal complement $\H^n$.
Writing $B$ in block form as $ \left(\begin{array}{c|c}
A & v \\
\hline
 w& q
\end{array}\right)$ we therefore have $w=v=0$.
Consequently,
$
B=
 \left(\begin{array}{c|c}
A & 0 \\
\hline
0 & q
\end{array}\right) 
$
with $
A\in Sp(n), q\in Sp(1),
$
which proves \eqref{N1}.\negro

The dimensions are: 
\begin{align*}
\dim G =& (n+1)(2n+3) & \dim H=&n(2n+1) &  \dim N =& n(2n+1) +3\\   G/H =& S^{4n+3} &  W=&S^3&
&&
\end{align*}

A generator of the Weyl group $W$ is induced by  an element of the form $\left(\begin{array}{c|c}
I & 0 \\
\hline
0 & q
\end{array}\right)$ with $q\in Sp(1)$. Since a diffeomorphism between $G/H$ and $S^{4n+3}$ is given by $AH \mapsto A  e_{n+1}$, the generator above acts on the sphere $S^{4n+3}$ by multiplication by $q$ on all the quaternionic coordinates. Then 
\begin{equation}\label{HP}
G/N = (G/H)/W = \H P^n.
\end{equation}

 Recall that the Lie algebra of $G$ is $\g =  \{ A \in \mathcal M_{n+1}(\H) \mid A^\dag= - A\}.$
Consider the $\mathrm{Ad}$-invariant inner product on $\mathfrak{g}$ defined by
$
\langle u, v \rangle = -\mathrm{Re}\,\mathrm{tr}(uv)$ for $u, v \in \mathfrak{g}.
$
This induces a Riemannian metric $\nu$ on $G$.

\medskip

The Lie algebra $\mathfrak{g}$ admits the orthonormal  basis
\[
\mathcal{B} = \{u_1, \ldots, u_m,\; v_1, v_2, v_3,\; t_{1,1}, \ldots, t_{1,n},\; t_{2,1}, \ldots, t_{2,n},\; t_{3,1}, \ldots, t_{3,n},\; t_{4,1}, \ldots, t_{4,n}\},
\]
with $m=n(2n+1)$, where
\begin{itemize}
\item[+] $u_1 = 
\left(\begin{array}{c|c}
u_1 & 0 \\
\hline
0 & 0
\end{array}\right), 
\ldots , u_{m} = \left(\begin{array}{c|c}
u_{m} & 0 \\
\hline
0 & 0
\end{array}\right)$
where $\{u_1, \ldots, u_{m}\}$ is any basis of $\h$.

\item[+] $v_1  =
\frac 1 {2} \left(\begin{array}{c|c}
0 & 0 \\
\hline
0 & i
\end{array}\right), v_2  =
\frac 1 {2} \left(\begin{array}{c|c}
0 & 0 \\
\hline
0 & j
\end{array}\right), 
v_3  =
\frac 1 {2} \left(\begin{array}{c|c}
0 & 0 \\
\hline
0 & k
\end{array}\right), 
$

\item[+] $t_{1,1} =  \frac 1 2 \left(\begin{array}{c|c}
0 & e_1 \\
\hline
-e_1^* & 0
\end{array}\right),
\ldots,
t_{1,n} =  \frac 1 2 \left(\begin{array}{c|c}
0 & e_n \\
\hline
-e_n^* & 0
\end{array}\right).
$

\item[+] $t_{2,1} =  \frac 1 2i \left(\begin{array}{c|c}
0 & e_1 \\
\hline
 e_1^* & 0
\end{array}\right),
\ldots,
t_{2,n} =  \frac 1 2i \left(\begin{array}{c|c}
0 & e_n \\
\hline
e_n^* & 0
\end{array}\right).
$

\item[+] $t_{3,1} =  \frac 1 2j \left(\begin{array}{c|c}
0 & e_1 \\
\hline
 e_1^* & 0
\end{array}\right),
\ldots,
t_{3,n} =   \frac 1 2 j \left(\begin{array}{c|c}
0 & e_n \\
\hline
e_n^* & 0
\end{array}\right).
$
\item[+] $t_{4,1} =  \frac 1 2k \left(\begin{array}{c|c}
0 & e_1 \\
\hline
 e_1^* & 0
\end{array}\right),
\ldots,
t_{4,n} =   \frac 1 2k\left(\begin{array}{c|c}
0 & e_n \\
\hline
e_n^* & 0
\end{array}\right).
$
\end{itemize}
Here, $\{e_1, \ldots, e_n\}$ denotes the  canonical basis of $\H^n$ and $e^*_j$ is the conjugate transpose of $e_j$. 
We use the decomposition
$$
\g = \h \oplus \z \oplus \m,
$$
where $\h =\mathrm{span}\{u_1, \ldots, u_m\}, \z =\mathrm{span}\{ v_1, v_2,v_3\}$ and
$\m = \mathrm{span}\{ t_{1,1} \ldots t_{4,n}\}$. Notices that 
$\n =\h \oplus {\z}$ is the normalizer of $\h$, that is, the Lie algebra of $N$.
The brackets between the elements of this decomposition are described as follows.

\begin{enumerate}[label=(\roman*)]
\item $[\h,\h] \subset \h$,
\item $[\h,\z] = [\z,\h] =0$, 
\item $[\h,\m] = [\m,\h] \subset \m$, 
\item $[v_1,v_2] = v_3, [v_1,v_3] = - v_2, [v_2,v_3] = v_1$,
\item  $[\z,\m] =[\m,\z] \subset \m$ is given by
\[
\begin{array}{c|cccc}
 [- ,- ] &  t_{1,\ell} & t_{2,\ell}  & t_{3,\ell} & t_{4,\ell}  \\\hline
v_1&  -t_{2,\ell} &  t_{1,\ell} &t_{4,\ell} &- t_{3,\ell} \\
v_2&  -t_{3,\ell}& -t_{4,\ell}&t_{1,\ell}  &t_{2,\ell} \\
 v_3 & -t_{4,\ell} &t_{3,\ell}&-t_{2,\ell} &t_{1,\ell} \\
\end{array}
\]
multiplied by $\frac 1 {2}$, where $\ell \in \{1, \ldots,n\}$,

\item the brackets of basis elements of $\mathfrak{m}$ lie in $\mathfrak{h}$, except for the following cases, 
for which the $\z$-component is given by $y$:
\[
\begin{array}{c|cccc}
y &  t_{1,\ell} & t_{2,\ell}  & t_{3,\ell} & t_{4,\ell}  \\
\hline
t_{1,\ell} &  &  -v_1 &-v_2& -v_3 \\
t_{2,\ell}&  v_1& &v_3 &-v_2\\
 t_{3,\ell} & v_2&-v_3&&v_1 \\
 t_{4,\ell} & v_3 &v_2&-v_1& \\
\end{array}
\]
where $\ell \in \{1, \ldots,n\}$.
\end{enumerate}

\begin{parrafo}{\bf Fundamental vector fields and dual differential forms }\label{61}.
The family of left invariant vectors fields associated to the elements of  $\mathcal B$ is 
$$Y_1, \ldots, Y_{m-1}, Q_1, Q_2 , Q_3,T_{1,1},\ldots,T_{1,n} , T_{2,1},\ldots,T_{2,n} , T_{3,1},\ldots,T_{3,n} , T_{4,1},\ldots,T_{4,n}.
$$
 The dual $\nu$-differential forms are denoted respectively by
$$
\gamma_1, \ldots \gamma_{m-1}, \lambda_1,\lambda_2, \lambda_3, \tau_{1,1}, \ldots , \tau_{1,n}, \tau_{2,1}, \ldots , \tau_{2,n}, \tau_{3,1}, \ldots , \tau_{3,n},\tau_{4,1}, \ldots , \tau_{4,n}.
$$
A straightforward calculation, using the previous bracket properties,  gives 
\begin{enumerate}[label=(\alph*)]
\item The Lie derivatives $L_{Y_\bullet} \lambda_\bullet$ are 0.

\item The Lie derivatives $L_{Q_\bullet} \lambda_\bullet$ are given by\[
\begin{array}{c|ccc} 
 L_{Q_\bullet} \lambda_\bullet&  \lambda_1 & \lambda_2 & \lambda_3  \\\hline
Q_1& 0 & - \lambda_3 &  \lambda_2 \\
Q_2&   \lambda_3&0 &- \lambda_1 \\
Q_3&- \lambda_2& \lambda_1& 0. 
\end{array}
\]
%
%
%
%
\item the differentials associated to $\lambda_\bullet$ are
\begin{align*}
    d\lambda_1 &= \lambda_2 \wedge \lambda_3 + \kappa_1, &
    d\lambda_2 &= -\lambda_1 \wedge \lambda_3 + \kappa_2, &
    d\lambda_3 &= \lambda_1 \wedge \lambda_2 + \kappa_3, \nonumber \\
    d\kappa_1 &= -\kappa_2\wedge\lambda_3+\kappa_3\wedge\lambda_2, &
    d\kappa_2 &= \kappa_1\wedge\lambda_3-\kappa_3\wedge\lambda_1, &
    d\kappa_3 &= -\kappa_1\wedge\lambda_2+\kappa_2\wedge\lambda_1 \label{gammakappa},
\end{align*}
\hspace{-.5cm} where
$\displaystyle \kappa_1 = \sum_{\ell=1}^n  (-\tau_{1,\ell} \wedge  \tau_{2,\ell}  +   \tau_{3,\ell}  \wedge \tau_{4,\ell} ),
\kappa_2 = 
\sum_{\ell=1}^n( -\tau_{1,\ell} \wedge  \tau_{3,\ell}  -   \tau_{2,\ell}  \wedge \tau_{4,\ell} ) $ and \\
$\displaystyle 
\kappa_3 =  \sum_{\ell=1}^n (-\tau_{1,\ell}  \wedge  \tau_{4_\ell} +    \tau_{2,\ell}  \wedge \tau_{3,\ell}  ).
$

 \item The Lie derivatives $L_{Q_\bullet} \kappa_\bullet$ are 
 \[
\begin{array}{c|ccc}
 L_{Q_\bullet} \kappa_\bullet&  \kappa_1 & \kappa_2 & \kappa_3  \\
\hline
Q_1& 0 &  -\kappa_3 & \kappa_2 \\
Q_2&   \kappa_3&0 & -\kappa_1 \\
Q_3&  -\kappa_2& \kappa_1& 0 \\
\end{array}
\]
\end{enumerate}
\end{parrafo}


\begin{parrafo}{\bf The $W$-action  on $G/H$.}\label{52} 
The Weyl group \(W=S^3\) acts on \(G/H=S^{4n+3}\) on the right, with quotient \(G/N=\mathbb{H}P^n\) (cf. \eqref{HP}). Recall that the cohomology of \(G/N\) can be computed using the complex of \(G\)-left-invariant differential forms:
\begin{equation}\label{GN}
\coho \BOm *{G/N}
=
\left\{
\omega\in
\bigwedge\nolimits^*(\tau_{1,1},\ldots,\tau_{4,n})
\ \middle|
\substack{
L_{Y_j}\omega=0,\quad 1\leq j\leq m,\\
L_{Q_k}\omega=0,\quad 1\leq k\leq 3
}
\right\}.
\end{equation}
The fundamental vector fields associated with the action of $W$ are $\{Q_1,Q_2,Q_3\}$ and their corresponding dual differential forms are $\lambda_1,\lambda_2,\lambda_3$. We also consider the differential forms $\kappa_1, \kappa_2, \kappa_3\in \lau \Om4{hor}{G/H}  $. The cycle
$
\kappa= \kappa_1^2 + \kappa_2^2 + \kappa_3^2 \in \coho \BOm 4{G/N}
$
(cf. \eqref{GN}, (a), (c) and (d)) defines the {\em Euler class} $[\kappa] \in \coho H4{G/N}$. 
It induces the map
$\delta\colon\coho H*{G/N}\to\coho H{*+4}{G/N}$
 given by multiplication by $[\kappa]$. This map is the connecting morphism of the Gysin sequence associated to the Weyl group action on $G/H$ and therefore $[\kappa]\ne 0$.

This class generates the cohomology of this quaternionic projective space
$G/N$. Thus, the volume form is a multiple of $\kappa^n$. Hence, there exists $C_n\in \R\menos \{0\}$ such that 
\begin{equation}\label{123}
 \kappa_{vol} =  \tau_{1,1} \wedge \tau_{2,1} \wedge \tau_{3,1} \wedge \tau_{4,1}  \wedge \cdots \wedge  \tau_{1,n} \wedge \tau_{2,n} \wedge \tau_{3,n} \wedge \tau_{4,n} = C_n \kappa^n.
\end{equation}
\end{parrafo}

The cycle $\kappa$ defines a cohomology class in $\coho H 4 M$ 
(cf. \eqref{FormulaHiru}), which we study in \lemref{kappae}.

\begin{parrafo}{\bf The $W$-action on $M^H$.}\label{kappatau}
\label{Euler3} 
The Weyl group $W=S^3$ acts almost freely on the manifold $M^H$ and, as we will see, the associated Euler form plays a key role in the Gysin sequence constructed in this work.

  We denote by $\{X_1,X_2,X_3\}$ the fundamental vector fields of the action and by $\{\chi_1, \chi_2,\chi_3\}$ their dual forms. The derivatives arise in this way,
\begin{align}\label{xe}
d\chi_1 &= e_1+\chi_2\wedge\chi_3, 
& d\chi_2 &= e_2-\chi_1\wedge\chi_3, 
& d\chi_3 &= e_3+\chi_1\wedge\chi_2,\\
de_1 &= -e_2\wedge\chi_3+e_3\wedge\chi_2, 
& de_2 &= e_1\wedge\chi_3-e_3\wedge\chi_1, 
& de_3 &= -e_1\wedge\chi_2+e_2\wedge\chi_1. \nonumber
\end{align}
where $e_1,e_2,e_3 \in \lau\Om *{hor} E$.

The Euler form  $e = e_1^2+e_2^2+e_3^2 \in \coho \Omega4{E/N}=  \coho \Omega4{M/G}$ is a  basic cycle defining the {\em Euler class} $\epsilon \in \coho H4{E/N} =  \coho H4{M/G}$.
It induces the map
$\delta\colon\coho H*{M/G}\to\coho H{*+4}{M/G}$
given by multiplication by $[e]$. This map is the connecting morphism of the Gysin sequence associated to the Weyl group action.

Throughout the sequel, we will make use of the following formul\ae\ involving Lie derivatives.

\begin{equation}\label{22BC}
\begin{array}{c@{\qquad\qquad}c}
\begin{array}{c|ccc}
L_{X_i}\chi_\bullet & \chi_1 & \chi_2 & \chi_3 \\ \hline
X_1 & 0 & -\chi_3 & \chi_2 \\
X_2 & \chi_3 & 0 & -\chi_1 \\
X_3 & -\chi_2 & \chi_1 & 0
\end{array}
&
\begin{array}{c|ccc}
L_{X_i}e_\bullet & e_1 & e_2 & e_3 \\ \hline
X_1 & 0 & -e_3 & e_2 \\
X_2 & e_3 & 0 & -e_1 \\
X_3 & -e_2 & e_1 & 0
\end{array}
\end{array}
\end{equation}

\end{parrafo}

The classes $[\kappa], [e] \in \coho H 4 M$ we have just defined are closely related.

\begin{lemma}\label{kappae}
The two classes $[\kappa]$ and $[e]$ are equal in $\coho H4M$.
\end{lemma}
\begin{proof}
Let us consider the differential form 
\begin{align*}\label{theta}
\theta = &( \lambda_1-\chi_1)\wedge(\lambda_2-\chi_2)\wedge(\lambda_3-\chi_3)+(\lambda_1-\chi_1)\wedge (\kappa_1+e_1) \\ &+(\lambda_2-\chi_2)\wedge(\kappa_2+e_2) +(\lambda_3-\chi_3)\wedge(\kappa_3+e_3). \nonumber
\end{align*}
It verifies
\begin{itemize}
\item $\theta \in \coho \BOm 3 M$, and
\item $\kappa - e = d\theta$.
\end{itemize}
For the first property it suffices to apply \eqref{FormulaHiru}, (a), (b), (c), (d) and \eqref{22BC}. For the second property we use the notation
$$
\beta_\ell = \lambda_\ell -\chi_\ell,
$$
$\ell=1,2,3$, and we write
\begin{align*}
\theta_1& = 
\beta_1 \wedge\beta_2 \wedge \beta_3 \\
\theta_2 &= \beta_1 \wedge (\kappa_1+e_1) +  \beta_2 \wedge (\kappa_2+e_2) + \beta_3 \wedge (\kappa_3+e_3) ,
\end{align*}
which gives $\theta =\theta_1+\theta_2$.
To obtain $d\theta_1$, we need the following calculation.
\begin{align*}
d \beta_1 \wedge \beta_2\wedge \beta_3
 & =_{(c),\eqref{xe}}  (\lambda_2\wedge \lambda_3 -\chi_2\wedge \chi_3 + \kappa_1-e_1)\wedge(\lambda_2 -\chi_2) \wedge (\lambda_3 -\chi_3) = (\kappa_1 - e_1)  \wedge \beta_2 \wedge \beta_3,  \\
 \beta_1 \wedge  d\beta_2\wedge \beta_3&=  -(\kappa_2- e_2)  \wedge\beta_1\wedge \beta_3, \\
\beta_1 \wedge \beta_2\wedge d\beta_3& = (\kappa_3 - e_3)  \wedge \beta_2 \wedge \beta_3.
\end{align*}
We get the first derivative
\begin{align*}
d\theta_1&= 
(\kappa_1 - e_1 ) \wedge \beta_2 \wedge \beta_3 + (\kappa_2 -e_2)  \wedge\beta_1 \wedge \beta_3  +(\kappa_3 - e_3 ) \wedge\beta_1 \wedge \beta_2.
 \end{align*}
To obtain $d\theta_2$, we need the following calculations.
\begin{align*}
d \beta_1 \wedge  (\kappa_1+e_1)  -  \beta_1 \wedge d (\kappa_1+e_1)  &=_{(c),\eqref{xe}}
(\lambda_2\wedge \lambda_3 -\chi_2\wedge \chi_3)\wedge(\kappa_1+e_1)
+
\kappa^2_1 -e^2_1\\
&-
\beta_1 \wedge (- \kappa_2\wedge \lambda_3+\kappa_3\wedge \lambda_2 - e_2 \wedge \chi_3 +e_2\wedge \chi_1 ), \\
d \beta_2 \wedge  (\kappa_2+e_2)  -  \beta_2 \wedge d (\kappa_2+e_2)  &=_{(c),\eqref{xe}}
-(\lambda_1\wedge \lambda_3 -\chi_1\wedge \chi_3)\wedge(\kappa_2+e_2)
+
\kappa^2_2 -e^2_2 \\
&+
\beta_2 \wedge (- \kappa_1\wedge \lambda_3+\kappa_3\wedge \lambda_1 - e_1 \wedge \chi_3 +e_1\wedge \chi_1 ) ,\\ 
d \beta_3 \wedge  (\kappa_3+e_3)  -  \beta_3 \wedge d (\kappa_3+e_3)  &=_{(c),\eqref{xe}}
(\lambda_1\wedge \lambda_2 -\chi_1\wedge \chi_2)\wedge(\kappa_3+e_3)
+
\kappa^2_3 -e^2_3 \\
&-
\beta_3 \wedge (- \kappa_1\wedge \lambda_2+\kappa_2\wedge \lambda_1 - e_1 \wedge \chi_2 +e_2\wedge \chi_1 ). 
\end{align*}
And then we get the second derivative
 \begin{align*}
d\theta_2 &=
( -\chi_2\wedge \chi_3  +\chi_2\wedge \lambda_3 + \lambda_3\wedge \lambda_2-\chi_3\wedge\lambda_2)\wedge (\kappa_1 - e_1) \\
&+ (- \chi_1\wedge \lambda_3 +\chi_1\wedge \chi_3 -\lambda_3\wedge \lambda_1 + \chi_3\wedge \lambda_1)\wedge (\kappa_2 - e_2)\\
& + (-\lambda_1\wedge \lambda_2 + \chi_1\wedge \lambda_2  -\chi_2\wedge \lambda_1 -\chi_1\wedge \chi_2)\wedge (\kappa_3 - e_3) + \kappa-e
\\
 &=
 -\beta_2 \wedge \beta_3 \wedge (\kappa_1-e_1)  -\beta_1 \wedge \beta_3\wedge (\kappa_2-e_2) -\beta_1 \wedge \beta_2\wedge (\kappa_3-e_3)  + \kappa-e.
  \end{align*}
 Finally,
$$
d\theta = d\theta_1 + d \theta_2 = \kappa -e.
\qedhere
$$
\end{proof}

\bigskip

\centerline {
\fbox{\begin{tabular}{c}
Background on Case II\\
$(G,H) = (SU(n+1), SU(n))$, with $n\geq 2$, \text{or} \\
$
(G,H) =  (U(n+1) , U(n))
\quad\text{or}\quad
(G,H) = (Sp(n+1)\cdot U(1), Sp(n) \cdot U(1))$, with $n\geq 1$.
\end{tabular}}
}

\bigskip

Let us begin with the Lie groups. We have the following properties.

\bigskip

\noindent {\bf Case II-a}. $(G,H) = (SU(n+1), SU(n))$ with $n\geq2$. 
Recall that $G=  \{ A \in \mathcal M_{n+1}(\C) \mid A^*A=I \hbox{ and } \det A =1\}$. The inclusion $H \subset G$ identifies $A\in H$ with 
$\left(\begin{array}{c|c}
A & 0 \\
\hline
0 & 1
\end{array}\right)\in G$. 
We claim that, in fact,
\begin{equation*}\label{N2a}
N = \left\{ \left(\begin{array}{c|c}
A & 0 \\
\hline
0 & (\det A)^{-1}
\end{array}\right)
\mid A \in U(n)\right\} \cong U(n).
\end{equation*}
The proof is analogous to  the one found in the Appendix - Background on Case III.
We have
\begin{align*}
\dim G =& (n+1)^2 -1,& \dim H=&n^2 -1,&  \dim N =& n^2, \\   G/H =& S^{2n+1}, &  W=&S^1.&
&&
\end{align*}

\medskip

\noindent {\bf Case II-b}.  We have $(G,H) = (U(n+1), U(n))$ with $n\geq1$. 
Recall that $G=  \{ A \in \mathcal M_{n+1}(\C) \mid A^*A=I\}$. The inclusion $H \subset G$ identifies $A\in H$ with 
$\left(\begin{array}{c|c}
A & 0 \\
\hline
0 & 1
\end{array}\right)\in G$. 
We claim that, in fact, the normalizer $N$ of $H$ in $G$ is
$$
N = \left\{ \left(\begin{array}{c|c}
A & 0 \\
\hline
0 & \lambda
\end{array}\right)
\middle| A \in H, \lambda \in U(1)\right\} \cong H\times U(1).
$$
The proof is analogous to  the one found in the Appendix - Background on Case III.
We have
\begin{align*}
\dim G =& (n+1)^2, & \dim H=&n^2, &  \dim N =& n^2 +1,\\   G/H =& S^{2n+1}, &  W=&S^1.&
&&
\end{align*}

\medskip

\noindent {\bf Case II-c}. We have $(G,H) = (Sp(n+1)\cdot U(1), Sp(n)\cdot U(1))$ with $n\geq1$.
Recall that $G= \dfrac{Sp(n+1) \times U(1)}{\{(I,1), (-I,-1)  \}}$. The inclusion $H \subset G$ identifies the class $[A,z] \in H$ with the class
$\left[\left(\begin{array}{c|c}
A & 0 \\
\hline
0 & z
\end{array}\right),z\right]\in G$. 
We claim that, in fact, the normalizer $N$ of $H$ in $G$ is
\begin{equation}\label{N2}
N = \left\{ \left[ \left(\begin{array}{c|c}
A & 0 \\
\hline
0 & \lambda
\end{array}\right),\mu \right]
\middle|  [A,\mu] \in H, \lambda \in U(1) \right\} \cong Sp(n) \cdot U(1)  \cdot U(1).
\end{equation}
Consider an  element $[C,q] \in N$. 
Since $N$ normalizes $H$, it normalizes its $Sp(n)$-factor. Hence, $C \in N_{Sp(n+1)}(Sp(n))$. 
Therefore, 
$C
=
\left( \begin{array}{c|c}
B & 0 \\
\hline
0 & p
\end{array}\right)$
for some $B \in Sp(n)$ and $p\in Sp(1)$ (see Appendix - Background Case III).
The condition  $[C,q]\in N$ is tantamount to
$$
 \left[\left(\begin{array}{c|c}
B & 0 \\
\hline
0 & p
\end{array}\right), q\right]
\underbrace{ \left[\left(\begin{array}{c|c}
A & 0 \\
\hline
0 &z
\end{array}\right),z\right]}_\alpha
 \left[\left(\begin{array}{c|c}
B^{-1} & 0 \\
\hline
0 & p^{-1}
\end{array}\right), q^{-1}\right]
\in H
$$
for each $\alpha \in H$. As $pzp^{-1}=qzq^{-1}$ for every $z\in U(1)$, assertion $\eqref{N2}$ follows.

We have
\begin{align*}
\dim G =& (n+1)(2n+3) +1 & \dim H=&n(2n+1) +1&  \dim N =& n(2n+1) +2\\   G/H =& S^{4n+3} & W=&S^1&
&&
\end{align*}

The following result summarizes the common properties of the Case II actions.

\begin{proposition} Consider a regular spherical action of $G$ on a manifold $M$ in Case II. 
Then the following properties hold.

(II-1) The Weyl group $N/H$ is a circle,

(II-2) there exists a $Ad$-invariant scalar product $\langle -,-\rangle$ in $\g$, and

(II-3) there exists a orthonormal basis $\mathcal B = \{u_1, \ldots , u_{m},u, v_1, \ldots, v_p,w_1, \ldots , w_p\}$ of $\g$ with  $\{u_1, \ldots , u_{m}\}$ basis of $\h$ and  $\{u_1, \ldots , u_{m},u\}$ basis of $\n$.

We have the following decomposition of $\g$ as a direct sum of vector spaces:
$$
\g =  \underbrace{\langle u_1, \ldots, u_{m} \rangle}_{\h} \oplus  \underbrace{\langle u \rangle}_{\z}\oplus  \underbrace{\langle  v_1, \ldots, v_p,w_1, \ldots , w_p \rangle}_{\m}.
$$

(III-4) There exist constants $K_1, \ldots, K_p \in \R$, $C_1, \ldots C_p\in \R\menos \{ 0 \}$, with

\begin{enumerate}[label=(\roman*)]
\item $[\h,\h] \subset \h$,
\item $[\h,\z] = [\z,\h] =0$, 
\item $[\h,\m] = [\m,\h] \subset \m$, 
\item $[v_\ell,u ] =  K_\ell w_\ell \hbox{ and } [w_\ell,u ] =   - K_\ell v_\ell $, for each $\ell \in \{1,\ldots,p\}$,
\item$[v_\ell,w_\ell] \in  \h + C_\ell  \  \z$, for each $\ell\in \{1, \ldots, p\}$, and
\item $\{ [v_l,v_\ell], [v_l,w_\ell], [w_l,v_\ell], [w_l,w_\ell] \}  \subset \h$, for each $l,\ell\in \{1,\ldots,p\}$ except the previous case.
\end{enumerate}
\end{proposition}

\begin{proof} We consider each of the three families in Case III separately.

\medskip

{\bf Special unitary group}.  \label{SU} The first example of a generic couple  arrives when $(G,H) = (SU(n+1), SU(n))$ with $n\geq 2$, that is, Case II-a.

Let us verify properties (II-1) - (II-4). Following \tabref{devito} we have (II-1).
We set $m=n^2-1,p= n$.

Recall that the Lie algebra of $G$ is $\g = \{ A \in \mathcal M_{n+1}(\C) \mid A^* =-A \hbox{ and } \tr A =0\}$
 We consider on $\g$ the scalar product
$
\langle u,v\rangle = -\Re \  \tr (uv) ,
$
where
$u,v\in \g$, 
which is invariant under adjunction. This gives (II-2).
The basis $\mathcal B = \{ u_1, \ldots, u_m, u,v_1, \ldots v_p,w_1, \ldots, w_p\}$ we are looking for is
\begin{itemize}
\item $u_1 = 
\left(\begin{array}{c|c}
t_1 & 0 \\
\hline
0 & 0
\end{array}\right), 
\ldots , u_{m} = \left(\begin{array}{c|c}
t_{m} & 0 \\
\hline
0 & 0
\end{array}\right)$
where $\{t_1, \ldots, t_{m}\}$ is any orthonormal basis of $\h$,

\item $u = \frac 1 {\sqrt{n^2+n}} i \left(\begin{array}{c|c}
I & 0 \\
\hline
0 & -n
\end{array}\right)$,
 where $I\in \h$ is the identity matrix,

\item$v_1 =\frac 1{ \sqrt 2}\left(\begin{array}{c|c}
0 & e_1 \\
\hline
-e_1^* & 0
\end{array}\right),
\ldots,
v_n = \frac 1{ \sqrt 2}\left(\begin{array}{c|c}
0 & e_n \\
\hline
-e_n^* & 0
\end{array}\right),
$

\item $w_1 = \frac 1{ \sqrt 2} i \left(\begin{array}{c|c}
0 & e_1 \\
\hline
e_1^* & 0
\end{array}\right),
\ldots,
w_n =\frac 1{ \sqrt 2} i\left(\begin{array}{c|c}
0 & e_n \\
\hline
e_n^* & 0
\end{array}\right).
$
\end{itemize}
Here, $\{e_1, \ldots, e_n\}$ denotes the  canonical basis of $\C^n$.  A straightforward calculation gives that properties (II-3) and (II-4) are verified with 
$ K_\ell = -C_\ell= - \sqrt {\frac{n+1} {n}}$ and  $\ell\in\{1,\ldots,n\}$.

\medskip 

{\bf Unitary group}.  \label{UU} The second example of a generic couple arrives when $(G,H) = (U(n+1), U(n))$ with $n\geq1$, that is,  Case II-b.

Let us verify properties (II-1) - (II-4). Following \tabref{devito} we have (II-1).
We set $m=n^2, p= n$.

Recall that the Lie algebra of $G$ is $\g =  \{ A \in \mathcal M_{n+1}(\C) \mid A^*= - A\}$. We consider on $\g$ the scalar product
$
\langle u,v\rangle = - \Re \ \tr (uv) ,
$
where $u,v\in \g$, 
which is invariant under adjunction. 
This gives (II-2).
The basis $\mathcal B = \{ u_1, \ldots, u_m, u,v_1, \ldots v_p,w_1, \ldots, w_p\}$ we are looking for is
\begin{itemize}
\item$u_1 = 
\left(\begin{array}{c|c}
t_1 & 0 \\
\hline
0 & 0
\end{array}\right), 
\ldots , u_{m} = \left(\begin{array}{c|c}
t_{m} & 0 \\
\hline
0 & 0
\end{array}\right)$
where $\{t_1, \ldots, t_{m}\}$ is any orthonormal basis of $\h$,

\item $u = \frac 1 {\sqrt {n+1}} i I$, where $I\in \g$ is the identity matrix,

\item $v_1 = \frac 1 {\sqrt 2} \left(\begin{array}{c|c}
0 & e_1 \\
\hline
-e_1^* & 0
\end{array}\right),
\ldots,
v_n =  \frac 1 {\sqrt 2} \left(\begin{array}{c|c}
0 & e_n \\
\hline
-e_n^* & 0
\end{array}\right),
$

\item $w_1 =   \frac 1 {\sqrt 2}i \left(\begin{array}{c|c}
0 & e_1 \\
\hline
e_1^* & 0
\end{array}\right),
\ldots,
w_n = \frac 1 {\sqrt 2} i\left(\begin{array}{c|c}
0 & e_n \\
\hline
e_n^* & 0
\end{array}\right).
$
\end{itemize}
Here, $\{e_1, \ldots, e_n\}$ denotes the  canonical basis of $\C^n$. Notice that $u$ generates the center of the Lie algebra $\g$. A straightforward calculation gives that properties (II-3) and (II-4) are verified with 
$ K= 0$ and
$C_\ell=- \sqrt{n+1}$  with $\ell\in\{1,\ldots,n\}$.

\medskip

{\bf Central product}.  \label{UU} The third example of a generic couple appears when $(G,H) = (Sp(n+1)\cdot U(1), Sp(n)\cdot U(1))$ with $n\geq1$, that is,  Case II-c.

Let us verify properties (II-1), \ldots , (II-4). Following \tabref{devito} we have (II-1).
We set $m=n(2n+1)+1$ and $p= 2n+1$.

Recall that the Lie algebra of $G$ is $\g =\spp (n) \oplus \uu(1) =  \{ A \in \mathcal M_{n+1}(\H) \mid A^\dag =-A \} \oplus \R$.
 We consider on $\g$ the scalar product
$
\langle (u_1,u_2),(v_1,v_2) \rangle = -\Re \ \tr (u_1v_1)  -\Re \  \tr (u_2v_2) ,
$
with $(u_1,u_2),(v_1,v_2)\in \g$, 
which is invariant under adjunction. This gives (II-2).

The basis $\mathcal B = \{ u_1, \ldots, u_m, u,v_1, \ldots v_p,w_1, \ldots, w_p\}$ we are looking for is
\begin{itemize}
\item $u_1 = 
\left(\left(\begin{array}{c|c}
t_1 & 0 \\
\hline
0 & 0
\end{array}\right), 0\right)
\ldots ,
u_{m-1} = \left(\left(\begin{array}{c|c}
t_{m-1} & 0 \\
\hline
0 & 0
\end{array}\right), 0\right), u_{m} =    \left(0,1\right),
$
where $\{t_1, \ldots, t_{m-1}\}$ is any orthonormal basis of $\spp(n)$,

\item $u =  \left(\left(\begin{array}{c|c}
0& 0 \\
\hline
0 & i
\end{array}\right),0\right)$,

\item $v_1 = \left( \frac 1 {\sqrt{ 2}} \left(\begin{array}{c|c}
0 & e_1 \\
\hline
-e_1^* & 0
\end{array}\right),0\right),
\ldots,
v_n = \left( \frac 1 {\sqrt{ 2}} \left(\begin{array}{c|c}
0 & e_n \\
\hline
-e_n^* & 0
\end{array}\right),0\right),
$\\
$v_{n+1}  =\left(  \frac 1 {\sqrt{ 2}}  k\left(\begin{array}{c|c}
0 & e_1 \\
\hline
e_1^* & 0
\end{array}\right),0\right),
\ldots,
v_{2n} =  \left(  \frac 1 {\sqrt{ 2}}  k\left(\begin{array}{c|c}
0 & e_n \\
\hline
e_n^* & 0
\end{array}\right),0\right),
$
$
v_{2n+1} =  \left(   \left(\begin{array}{c|c}
0 & 0 \\
\hline
0 & k
\end{array}\right), 0\right),
$

\item $w_1 = \left( \frac 1 {\sqrt {2}}i \left(\begin{array}{c|c}
0 & e_1 \\
\hline
e_1^* & 0
\end{array}\right),0\right),
\ldots,
w_n = \left( \frac 1 {\sqrt{ 2}} i\left(\begin{array}{c|c}
0 & e_n \\
\hline
e_n^* & 0
\end{array}\right),0\right),
$\\
$w_{n+1}  =\left(  \frac 1 {\sqrt{ 2}}  j\left(\begin{array}{c|c}
0 & e_1 \\
\hline
e_1^* & 0
\end{array}\right),0\right),
\ldots,
w_{2n} =  \left(  \frac 1 {\sqrt{ 2}}  j\left(\begin{array}{c|c}
0 & e_n \\
\hline
e_n^* & 0
\end{array}\right),0\right),
$
$
w_{2n+1} =  \left(  \left(\begin{array}{c|c}
0 & 0 \\
\hline
0 & j
\end{array}\right), 0\right).
$

\end{itemize}
Here, $\{e_1, \ldots, e_n\}$ denotes the  canonical basis of $\H^n$.  A direct calculation gives that properties (II-3) and (II-4) are verified with 
$ K_\ell = -C_\ell =1$ and $K_{2n+1} = -C_{2n+1} =2$,  with $\ell\in\{1,\ldots,2n\}$.
\end{proof}
\begin{parrafo}
\label{55}{\bf Fundamental vector fields and dual differential forms} (see \parref{31}).
Consider a regular spherical action of $G$ on a manifold $M$ in Case II. 
The fundamental vector fields and the associated dual differential forms 
follow a common pattern, detailed below.

The family $\{Y_1, \ldots, Y_m, Q, V_1, \ldots , V_p,W_1,\ldots,W_p\}$ is the family of left invariant vectors fields of $G$ associated to the elements of  the basis $\mathcal B$. The family of dual differential forms of $G$, relatively to the scalar product given by property II-2,  is denoted by
$
\{\gamma_1, \ldots \gamma_m,\lambda,  \alpha_1, \ldots , \alpha_p, \beta_1, \ldots, \beta_p\}$.
A straightforward calculation, using the previous  properties (II-1)-(II-4),  gives 

\begin{enumerate}[label=(\Alph*)]
\item[]
\item[] $L_{Y_1} \lambda = \cdots = L_{Y_{m}} \lambda =L_Q\lambda = 0$. \hfill (A)
\item[]
\item[]  $L_{Q} \alpha_\ell=  K_\ell  \beta_\ell$ and $L_{Q} \beta_\ell = - K_\ell\alpha_\ell$ for $\ell\in \{1, \ldots, p\}$. \hfill (B)
\item[]
\item[] $L_{V_\bullet} \alpha_\bullet, L_{V_\bullet} \beta_\bullet, L_{W_\bullet} \alpha_\bullet, L_{W_\bullet} \beta_\bullet \in \mathrm{span}\{\gamma_1, \ldots,\gamma_m\}  $. \hfill (C)
\item[]
\item[]  $d\lambda =    \displaystyle\sum_{\ell=1}^{p} C_\ell\alpha_\ell \wedge \beta_\ell $ and 
$\left(d\lambda \right)^p = C \displaystyle \alpha_1 \wedge \beta_1\wedge \cdots \wedge \alpha_p \wedge \beta_p$, \hfill (D)
\end{enumerate}
where $C= C_1\cdots C_p\cdot  p!\ne 0$.

\end{parrafo}

\begin{parrafo}{\bf The $W$-actions on $M^H$ and $G/H$.}\label{Euler1} 
The Weyl group $W=S^1$ acts almost freely on the manifold $M^H$ and the associated Euler form plays a key role in the Gysin sequence constructed in this work.

 We denote by $X$ the fundamental vector field of the action and by $\chi$ its dual $1$-form. The differential form
$e=d\chi\in\coho\Om2{M^H/W}=\coho\Om 2{M/G}$
is a basic cycle since
\begin{equation}\label{Lchi}
L_X\chi=0.
\end{equation}
 The class
$\epsilon=[e]\in\coho H 2{M/G}$
is the \emph{Euler class of the Weyl-group action}. Multiplication by this class $\epsilon$ induces the connecting morphism 
$\delta\colon\coho H*{M/G}\to\coho H{*+2}{M/G}$ in the 
 Gysin sequence associated with the Weyl-group action on $M^H$ (see \cite[Proposition 14.33]{BottTu}).

The Weyl group $W=S^1$ acts on the homogeneous space $G/H$ by  right  multiplication
with $[d\lambda] \in \coho H 2 {G/N}$ as the Euler class of the associated Gysin sequence (see \cite[Proposition 14.33]{BottTu}).
These two corresponding classes induced on $M$ are the same, namely, $[d\lambda] = [e]$ on $\coho H2M$ since $d\lambda -e = d(\lambda-\chi)$ and $\lambda-\chi\in \coho \BOm 1 M$ (cf.   \eqref{FormulaBat}, (A) and \eqref{Lchi}).
\end{parrafo}

\centerline {
\fbox{\begin{tabular}{c}
Background on Case I\\

$(G,H) = (SO(n+1),SO(n))$, with $n\geq 2$ 
\quad\text{or}\quad
$(G,H) = (Sp(n+1) \cdot Sp(1), Sp(n)\cdot Sp(1))$, with $n\geq 1$ \quad\text{or}\quad
\\
$(G,H) = (G_2,SU(3))$ \quad\text{or}\quad $(G,H) = (Spin(7),G_2)$ \quad\text{or}\quad $(G,H) = (Spin(9), Spin(7))$
\end{tabular}}
}

\bigskip

In this case, we only need the following result.

\begin{proposition}\label{Z22}
In Case I, the Weyl group associated with the pair $(G,H)$ is $\mathbb{Z}_2$, and its generator acts on the sphere $G/H = S^d$ via the antipodal map.
\end{proposition}

\begin{proof}
The diffeomorphism between $G/H$ and $S^d$ of  \cite{MontSam} is constructed from the linear action $\Psi \colon G \times S^d \to S^d$ given by $A \mapsto  \Psi(A,x) =Ax$. We have $H=G_{x_0}$, where $x_0$ is a base point of the sphere,  and the diffeomorphism $f \colon G/H \to S^d$ defined by $f(AH) =Ax_0$. We proceed in three steps.

$\bullet$ {\em The orbit $W(x_0)$.} 
We have $f(N/H) = \left( S^d\right)^H$, as we now verify. If $B\in N$ and $C \in H$ there exists $C'\in H$ with  $CBx_0 = BC'x_0 = B x_0$ which proves the inclusion  $\subset$. Conversely, let $z \in S^d$ be such that $Cz =z$ for each $C\in H$. There exists $B\in G$ with $z=f(BH) = Bx_0$. Then, $Bx_0 = CBx_0 $ and therefore $B^{-1}CB \in G_{x_0} =H$ for each $C \in H$. It follows that $B\in N$ which proves the inclusion  $\supset$. The equality $f(N/H) = \left( S^d\right)^H$ implies that $f(N/H)$ is a sphere $S^0, S^1, \ldots$ containing $x_0$.

The Weyl group $W$ acts on $G/H$ from the right, inducing a free linear action $\Psi' \colon  W \times S^d \to S^d$  defined by $\Phi'(BH,x_0) = Bx_0$. By  \cite[Chapter~III, Proposition~8.5]{Bredon} we know that 
$W$ is isomorphic to either a finite group or one of the following three groups 
$S^1, N_{S^3}(S^1)$ or $S^3$. Since we have seen that the orbit $W(x_0) = f(N/H) $ is a sphere, it follows that  $W(x_0)$ is  one of $S^0,S^1$ or $S^3$ containing $x_0$.

$\bullet$ {\em The group $W$ is $\Z_2$.} We study each case.

{\bf Cases I-a, I-c and I-d}. In all three cases the inclusion $H\subset G$  is maximal among connected closed subgroups in the sense of inclusions (see for example \cite[Proposition~2.3]{GP}, \cite[Lemma 2.2.2]{LeMunir}, \cite[Lemma 2]{Varadarajan}). Therefore, $N_e=H$ or $G$. The latter case implies that $N=G$, meaning that, $H$ is a distinguished subgroup of $G$. Consequently, $S^d = G/H$ is a Lie group 
which is not possible
for $d=6,7$. Hence, $N_e = H$, and therefore $\dim W = 0$, which, by the previous item, implies that $W = \mathbb{Z}_2$.

\smallskip

\noindent {\bf Case I-c}. We have $(G,H) = (Spin(9), Spin(7)), n\geq 1$. Following \cite[Proposition~3.10]{Shankar}, we know that $N\subset Spin(8)$ . 
If we show that the inclusion $H \subset Spin(8)$ is maximal among connected closed subgroups, then, proceeding as in the previous step, we conclude that $W=\Z_2$. The maximality of the inclusion $H\subset Spin(8)$ can be established at the level of Lie algebras, namely, $\so(7) \subset \so(8)$. This is true by a result of Dynkin (see for example \cite[Proposition~2.3]{GP}).

\smallskip

\noindent {\bf Case I-e}. We have  $(G,H) = (Sp(n+1)\cdot Sp(1), Sp(n)\cdot Sp(1))$ with $n\geq1$.
Recall that $G= \dfrac{Sp(n+1) \times Sp(1)}{\Z_2}$. The inclusion $H \subset G$ identifies the class $[A,z] \in H$ with the class
$\left[\left(\begin{array}{c|c}
A & 0 \\
\hline
0 & z
\end{array}\right),z\right]\in G$. 
 We claim that, in fact,
\begin{equation}\label{NN}
N = \left\{ \left[ \left(\begin{array}{c|c}
A & 0 \\
\hline
0 & \lambda
\end{array}\right),\pm \lambda \right]
\mid [A,\pm\lambda] \in H, \lambda \in Sp(1) \right\} \cong Sp(n) \cdot Sp(1) \times \Z_2,
\end{equation}
which implies $W=\Z_2$. Let us verify \eqref{NN}.

Consider an  element $[C,\lambda] \in N$. 
Since $N$ normalizes $H$, it normalizes its $Sp(n)$-factor. Hence, $C \in N_{Sp(n+1)}(Sp(n))$. 
Therefore, 
$C
=
\left( \begin{array}{c|c}
B & 0 \\
\hline
0 & p
\end{array}\right)$
for some $B \in Sp(n)$ and $p\in Sp(1)$ (see Appendix - Background Case III).
The condition  $[C,\lambda]\in N$ is tantamount to
$$
 \left[\left(\begin{array}{c|c}
B & 0 \\
\hline
0 & p
\end{array}\right), \lambda\right]
\underbrace{ \left[\left(\begin{array}{c|c}
A & 0 \\
\hline
0 &z
\end{array}\right),z\right]}_\alpha
 \left[\left(\begin{array}{c|c}
B^{-1} & 0 \\
\hline
0 & p^{-1}
\end{array}\right), \lambda^{-1}\right]
\in H
$$
for each $\alpha \in H$. 
Equivalently,
$
\lambda^{-1}p\in Z(Sp(1))=\{\pm1\},
$
and hence
$p=\pm \lambda$. This gives \eqref{NN}.

$\bullet$ {\em Antipodal map.}  Since $W(x_0) =\Z_2$ is a sphere containing $x_0$ then  $W(x_0) = \{ x_0, -x_0\}$.
Let $J\in N$ be an element representing a generator of $W$. 
We consider the commutative diagram
$$
\xymatrix{
G/H \ar[r]^f \ar[d]^h& S^d \ar[d]^{h'} \\
G/H \ar[r]^f & S^d,
}
$$
where $h$ is the right multiplication by $JH$ and $h'$ is the induced map. 
We claim that $h'$ is the antipodal map. Indeed, for each $z\in S^d$, there exists $C \in G$ with $f(CH) = Cx_0 = z$. Therefore
$$
h'(z) = h'(f(CH)) = f(h(CH)) = f(CJ) = CJ  x_0 .
$$
If $Jx_0=x_0$ then $h'$, and hence $h$, is the identity. Since $J\not\in H$, this is impossible. Therefore, $Jx_0=-x_0$ which implies $h'=-\id$,
 the antipodal map on $S^d$.
\end{proof}

\begin{parrafo}{\bf Generic examples. Calculation}. 
We have seen in \parref{example} that any regular spherical action $\Phi\colon G\times M\to M$, with associated groups $H,N$ and $W$ can be reconstructed from a smooth manifold $E$ and an almost free action $\Phi_0\colon W\times E\to E$.
Once this has been done, we have already seen that $E=M^H$ and $E/W= M^H/W=M/G$.

The \teoref{nagusia} shows us how to compute the cohomology of the manifold $M$ using $\coho H* {M/G}$ and $\coho H * {M^H}$. In cases II and III, the connecting morphism is given by the Euler class $\epsilon \in \coho H {\dim W+1} { M/G}$. 
All these ingredients can be computed just using the manifold $E$ and the action $\Phi_0$.
 In fact, we have $\coho H * {M/G} = \coho H * {E/W}$ and $\coho H * {M^H} = \coho H * E$. 
 Moreover, the Euler class  $\epsilon$ is precisely the Euler class associated with the action $\Phi_0$.

For the calculation of the cohomology of $M$, we have that the formula \eqref{Z2} becomes
 \begin{equation}
\coho H *M = \coho H * {E/W} \oplus
\begin{cases}
\coho H {*-d}{E/W}, & \text{if } d \text{ is odd},\\[4pt]
\bigl(\coho H{*-d}{E}\bigr)^{-W}, & \text{if } d \text{ is even}.
\end{cases}
\end{equation}
 In Case II, the Gysin sequence \eqref{Gys2} becomes
 $$
 \ldots \to \coho H k {E/W} \to \coho Hk  M \to  \coho H {*-2p-1}{E/W}\xrightarrow{\delta}  \coho H {k+1} {E/W} \to \cdots ,
$$
where the connecting homomorphism $\delta$ is given by multiplication by  $\epsilon^{p+1}$,  with $\epsilon\in \coho H 2 {E/W}$ the Euler form of the action $\Phi_0$. 
Finally, in Case III, the Gysin sequence \eqref{Gys3} becomes
 $$
 \ldots \to \coho H k {E/W} \to \coho Hk  M \to  \coho H {*-4n-3}{E/W}\xrightarrow{\delta}  \coho H {k+1} {E/W} \to \cdots ,
$$
where the connecting homomorphism $\delta$ is given by multiplication by  $\epsilon^{n+1}$,  with $\epsilon\in \coho H 4 {E/W}$  the Euler form of the action $\Phi_0$. 

\end{parrafo}

\bibliographystyle{abbrv}
\bibliography{Ultimo}
  \end{document}